\documentclass[11pt]{amsart}

\usepackage{rotating}
\usepackage{amssymb}
\usepackage{amsxtra}
\usepackage{graphics}
\usepackage{latexsym}
\usepackage{xcolor}
\usepackage{amsmath}
\usepackage{amssymb,amsthm,amsfonts}
\usepackage{mathtools}
\usepackage{amscd}
\usepackage[arrow, matrix, curve]{xy}
\usepackage{syntonly}
\usepackage{multirow}
\usepackage{fancyhdr}
\usepackage{hyperref}

\usepackage{etex}
\usepackage{tikz-cd}
\usetikzlibrary{matrix, calc, arrows}
\usepackage{braket}
\usepackage{adjustbox}
\usepackage[official]{eurosym}

\theoremstyle{plain}

\theoremstyle{definition}

\numberwithin{equation}{thm}

\newcommand{\rk}{{\rm rank}}

\usepackage{geometry}

\begin{document}\textcolor{white}{.}\\

\centerline{Hodge Theory and Higgs Bundle on Moduli Spaces of   Smooth Varieties } 
\centerline{ With Semi-Ample Canonical Line Bundles, Negativity and Hyperbolicity}\vspace{.2cm}
\centerline{Kang Zuo}\vspace{.2cm}
\section{\bf Introduction}.\\[.2cm]
This is an expository note based on my lecture in the conference Recent Developments in Hodge Theory
March 29-April 2, 2021 IMSA, UMIAMI.\\[.2cm]
Let $U$ be a complex quasi-projective manifold and $\bar Y=U\cup \bar S$ a smooth compactification,  we are particularly interested on $U$ as a base space $$f: V\to U$$ parametrizing smooth polarized varieties (with semi ample canonical line bundles). We detail here the notion the complex hyperbolicity of $U$ in its various aspects as well as the program motivated by the original Shafarevich conjecture.\\[.2cm] 
The powerful theory on variation of Hodge structures developed by P.A.  Griffiths (\cite{G1}, \cite{G2} and \cite{G3}) shows that the negativity of the sheaf of logarithmic holomorphic vector fields of the log pair $(\bar Y,\bar S)$ if the local Torelli injectivity holds.  C. Simpson  \cite{S1} has introduced a basic notion in nonabelian Hodge theory, the so-called graded Higgs bundle or the system of Hodge bundles $(E,\theta)$  as the grading of the Hodge filtration on a polarized variation of Hodge structures (PVHS).  This  $\mathcal{O} _Y$-linear object is better understandable in algebraic geometry and enjoy a crucial property in complex geometry. For example,  the kernel of the Higgs field is semi-negative.  
These theories are of fundamental importance in the study of the global geometry of families.\\[.2cm]
Given a non-isotrivial family $$f: X\to Y $$ of varieties with semi-ample canonical line bundles, we still hope the non-trivial variation of complex structures of the fibres will produce certain negativities of the log tangent sheaf on the base.\\ Assuming the reduced singular fibre $\Delta_{red}$ is a normal crossing diviosr over $S$. 
the classical logarithmic Kodaira-Spencer map
$$   \tau:T_Y(-\log S)\to R^1f_*T_{X/Y}(-\log \Delta)$$
introduced by Kodaira-Spencer
measures infinitesimal variations of complex structures on fibres parametrized by the log base $(Y,S)$. The works of Viehweg--Zuo \cite{VZ-1,VZ-2,VZ-3}, combined the above-mentioned theories by extending the map $\tau=:\tau^{n,0}$
to the higher direct images of wedge products of the relative log tangent sheaf
 $$\tau^{n-q,q}:T_Y(-\log S)\otimes R^{q}f_*T^q_{X/Y}(-\log\Delta)\xrightarrow{}R^{q+1}f_*T^{q+1}_{X/Y}(-\log\Delta),\quad 0\leq q\leq n$$ and introduced the so-called deformation Higgs bundle
$$(F,\tau):=( \bigoplus_q R^{q}f_*T^q_{X/Y}(-\log \Delta), \, \bigoplus_q\tau^{n-q,q}),$$
which is a logarithmetic graded Higgs bundle
over the log base $(Y,S)$.  Moreover, Viehweg--Zuo constructed a Higgs map 
$$\rho\colon (F,\tau)\to (E,\theta)\otimes A^{-1},$$ 
where $(E,\theta)$ is the Higgs bundle associated to the variation of Hodge structures on the middle cohomology of a new family  $g\colon Z\to Y$  built from a cyclic covering of $X$ by taking roots out of sections of the relative pluricanonical linear system twisted by an anti-ample line bundle $A^{-1}$ on $Y$.  Via the maximal non-zero iteration of the Kodaira--Spencer map of $\rho (F,\tau)$ (also known as the Griffiths--Yukawa coupling) we obtain a non-zero map
$$ \tau^m: S^mT_Y(-\log S)\xrightarrow{}
\text{ker}(\theta^{n-m-1,m+1}) \otimes A^{-1}$$ where $\text{ker}(\theta^{p,q})$ is the kernel of the Higgs field $\theta^{p, q}.$  As the kernel of the Higgs field on a graded Higgs bundle arising from a local system is semi-negative \cite{Z},
one obtains a generically ample subsheaf 
$$\mathcal A\hookrightarrow S^m\Omega^1_Y(\log S), $$
 where $\mathcal A$ is called generically ample
in the sense of Viehweg, i.e.
 for some symmetric power $S^n(\mathcal A)$ there exists
 a generically surjective map
 $$ A^{\oplus l }\to S^n(\mathcal A).$$
\\[.2cm]
In the analytic setting, one defines   a complex Finsler pseudo metric $ds^2_{\mathbb{ C} \,F}$ on $T_Y(-\log S)$ by taking the $m$-th root out of the product of the Hodge metric on the kernel of the Higgs field and the Fubini-Study metric on $A^{-1}$ via the above Griffiths-Yukawa coupling.
 One shows that  
$ds^2_{\mathbb{ C} \,F}$  has the holomorphic sectional curvature bounded above by the negative (1,1)-form of the Fubini-Study metric.  The complex Finsler metric with the negative holomorphic sectional curvature derived as above plays the same role as the negative holomorphic sectional curvature associated to horizontal period maps in Hodge theory.  It implies Brody and Kobayashi hyperbolities of $U$ in the works by \cite{VZ-1}, \cite{PTW} and \cite{Deng}. The generically ample subsheaf plays crucial role in the paper \cite{CP} on solving Viehweg conjecture. Also,   the big Picard theorem was recently proven in \cite{DLSZ} followed by applying Nevanlinna theorem to those type of complex Finsler pseudo metric. \\[.2cm]
Beside of the complex hyperbolicities the original Shafarevich program asked for the so-called finiteness of the set $\mathcal H$ of isomorphic classes of families over a fixed log base. The finiteness problem   is decomposed into two basic problems:
\begin{itemize}
	\item The boundedness of  $\mathcal H$.
	\item The rigidity of points in $\mathcal H$
\end{itemize}  
The boundedness has been proven by various people. For example, for semi-stable families of abelian varieties over a base curve Faltings \cite{F}  has 
shown the boundedness .  Jost-Yau  \cite{JY} have used Yau's form of Schwarz inequality and given another proof of the boundedness for families of abelian varieties over curves.  In general, Viehweg-Zuo \cite{VZ-2} have shown the boundedness for families of varieties with semi-ample canonical line bundles over a base curve without requiring the existence of a locally injective Torelli map.  Finally, Kov\'acs-Lieblich \cite{KL}   have shown the boundedness for families over a higher dimensional base.    \\
The rigidity problem is subtle. In fact, there exists non-rigid families of higher dimensional varieties.
So characterizing non-rigid families will be important for Shafarevich program for families of higher dimensional varieties. \\[.2cm]
The puroses of this note are\\[.2cm]
{\bf(1)} We review graded Higgs bundles
arising from polarized variation of Hodge structures and deformation (graded) Higgs bundles attached a logarithmic smooth family of varieties.\\[.2cm]
{\bf (2)} We give brief discussion on the constructions of the comparison map from the deformation Higgs bundle to the graded Higgs bundle from a polarized variation of Hodge structure and the big subsheaf in logarithmic symmetric differential forms on the log base $(Y,S).$ We also remark that the Viehweg-Zuo's construction extends to the case of families of varieties of general type or of good minimal model. In fact, it has been done already over 1-dimensional bases in \cite{VZ-1}.
\\[.2cm]
{\bf (3)} As the first application of the existence of the big subsheaf in logarithmic symmetric differential forms we give a uniform proof for the results on the algebraic hyperbolicity of $(Y,S)$ due to Parshin-Arakelov, Miglinorini, Zhang, 
Kov\'acs, Bedulev-Viehweg, Oguiso-Viehweg, and Viehweg-Zuo.\\[.2cm] 
{\bf (4)} We recall the classical Picard theorems, the so-called little Picard theorem and big Picard theorem in complex one variable and their generalizations   to higher dimensional varieties. We stetch the proof
 in the paper \cite{DLSZ} for
the Picard extension theorem on a base parametrizing smooth varierties with semi-ample canonical line bundles. \\[.2cm]
{\bf (5)} There are more applications of the comparison map. We discuss briefly the recent paper \cite{LYZ} on the strict Arakelov inequality of families of varieries of general type.\\[.2cm]
{\bf (6)} Very recently Javanpeykar-Sun-Zuo \cite{JSZ} proved a new version of the Shafarevich conjecture, the so-called finiteness of pointed families $f: X\to Y$ of $h$-dimensional polarized varieties. A log map from a log curve  $\phi: (C,S_C)\to (Y,S)$  has only the trivial deformation $\phi_t$ if $\phi_t$ fixes no less than $ {1\over 2} (h-1)\deg\Omega^1_C(\log S_C)$ number of points in $C\setminus S_C$.\\[.2cm]
{\bf 7) } We propose a program on characterization of non-rigid famlies.\\[.2cm]
{\bf 8)} We introduce the notion of toplogical hyperbolicity on a quasi-projective variety, and pose a conjecture on the topological hyperbolicity on moduli stacks of smooth projective varieties with semi-ample canonical line bundles.\\[.2cm]
\noindent{\bf Acknowledgment.}  I would like to thank Steven Lu   and Ruiran Sun for discussions on various notions of complex hyperboilicities, Riemann-Finsler metric and flag curvature.  I would also like to thank Jinbang Yang for disussions on Jacobian ring
of hypersurfaces in projective spaces and the relation to Hodge theory.

\section{ \bf Graded Higgs Bundle arising from Geometry}
\noindent
Throughout we will assume that $U$ is a quasi-projective manifold and compactified by a projective manifold $\bar Y$ with $\bar S=\bar Y\setminus U$ is a normal crossing divisor, and that there is smooth family $$f: V\to U$$ of $n$-folds. Leaving out a codimension two subset of $\bar Y$ we find a good partial compactification as a log smooth family $$f:X\to Y$$ in the following sense
\begin{itemize}
\item $X$ and $Y$ are quasi-projective manifolds, $f$ is flat, $U\subset Y$ and
$\mathrm{codim}(\bar Y\setminus Y)\geq 2.$
\item $S=Y\setminus U$ is smooth and $\Delta=f^*S$ is a relative normal crossing divisor over $S$ (i.e. whose reduced components, and all their intersections are smooth over components of $S$).
\end{itemize}.\\[.2cm]
\subsection{ \bf System of Hodge Bundles Attached to Polarized Variation of Hodge Structure}
Following
P.A.  Griffiths and C. Simpson one introduces the most natural Higgs bundle relating to geometry and topology on a log smooth family $$f: X\to Y.$$ 
The tautological sequence\\
{\bf (2.1.1)}
\[\quad 0\to f^*\Omega^1_Y(\log S)\to\Omega^1_X(\log\Delta)\to\Omega^1_{X/Y}(\log\Delta)\to 0\]\\
induces the short exact sequence of logarithmic forms of higher degrees\\
{\bf (2.1.2)}
\[ 0\to f^*(\Omega^1_Y)\otimes\Omega^{p-1}(\log S)\to\text{gr}\Omega^p_X(\log S)\rightarrow\Omega^p_{X/Y}(\log S)\to 0,\]
where
\[\text{gr}\Omega^p_{X}(\log\Delta)=\Omega^p_X(\log\Delta)/f^*\Omega^2_Y(\log S)\otimes\Omega^{p-2}_X(\log\Delta).\]\\
 By Simpson the direct sum of the direct image sheaves
$$E^{p,q}=R^qf_*\Omega^p_{X/Y}(\log\Delta),\quad p+q=k$$ endowed with the connecting maps in\\[.2cm] {\bf (2.1.3.)}
\[\theta^{p,q}: R^qf_*\Omega^p_{X/Y}(\log\Delta)
\xrightarrow{\partial}\Omega^1_Y(\log S)\otimes R^{q+1}f_*\Omega^{p-1}_{X/Y}(\log\Delta)\]\\
forms the so-called system of Hodge bundles 
\[(E,\theta)
=(\bigoplus_{p+q=k}E^{p,q},\bigoplus_{p+q=k}\theta^{p,q})\]
of weight-$k$  attached to the family.\\[.2cm]
Take the dual of {\bf (2.1.1)}
$$0\to T_{X/Y}(-\Delta)\to T_X(-\Delta)\to f^*T_Y(-\log S)\to 0$$
The connecting map of the direct image defines the logarithmic Kodaira-Spencer map
\[\tau: T_Y(-\log S)\to R^1f_*T_{X/Y}(-\log\Delta).\]
The Higgs field $\theta^{p,q}$ can be also defined as the cup-product with $\tau$
\[\theta^{p,q}: T_Y(-\log S)\otimes R^qf_*\Omega^p_{X/Y}(\log\Delta)\xrightarrow{\tau\otimes id}\]
\[R^1f_*T_{X/Y}(-\log\Delta)\otimes R^qf_*\Omega^p_{X/Y}(\log\Delta)\xrightarrow{\cup}
R^{q+1}f_*\Omega^{p-1}_{X/Y}(\log\Delta).\]
\\[.2cm]
{\bf Proposition 2.1.}{\sl
The Higgs bundle $(E,\theta)$ is the grading of Deligne's quasi-canonical extension of the variation of the polarized Hodge structure on $k$-th Betti cohomology $R^k_\text{B}f_*\mathbb Z_{V}$ of the smooth family $f:V\to U$.}
\\[.2cm]
Proposition 2.1. is well known and due to Griffiths using theory of harmonic forms in K\"ahler geometry.
The crucial point is the $E_1$-degeneration of Hodge to de Rham spectral sequence.
We present an algebraic approach due to Katz-Oda \cite{KO} It works   over any characteristic once we have $E_1$-degeneration of Hodge to de Rham spectral sequence.\\[.2cm]
The short exact sequence {\bf (2.1.2)} induces a short exact sequence of complexes of log differential forms\\[.2cm]
{\bf (2.1.4.)}
\[0\to f^*(\Omega^1_Y(\log S)\otimes\Omega^{\bullet}_{X/Y}(\log S)[-1]\to
K^0/K^2(\Omega^\bullet_{X/Y}(\log\Delta)\]
\[\to\Omega^{\bullet}_{X/Y}(\log S)\to 0.\]
The $k$-th direct image sheaf $V:=R_\text{dR}^k f_*\Omega^{\bullet}_{X/Y}(\log S)$ of the relative de Rham complex
carries the so-called logarithmic Gauss-Manin connection
\[V\xrightarrow{\nabla}V\otimes\Omega^1_Y(\log S)\]
defined by the connecting map in the long exact sequence
\\
\[\nabla: R_\text{dR}^k f_*\Omega^{\bullet}_{X/Y}(\log S)\to R_\text{dR}^{k+1}f_* f^*(\Omega^1_Y(\log S)\otimes\Omega^{\bullet}_{X/Y}(\log S)[-1]\]
\[\simeq\Omega^1_Y(\log S)
\otimes R^k_\text{dR}f_*\Omega^{\bullet}_{X/Y}(\log S).\]\\[.3cm]
Introducing a decreasing filtration of the truncated sub complexes
\[\{F_\text{tru}^p\Omega^\bullet_{X/Y}(\log\Delta)\}_{0\leq p\leq n}\subset\Omega^\bullet_{X/Y}(\log\Delta)\]
by
\[ 0\to\cdots\to 0\to\Omega^p_{X/Y}(\log\Delta)\xrightarrow{d}\Omega^{p+1}_{X/Y}(\log\Delta)\xrightarrow{d}\cdots\xrightarrow{d}\Omega^n_{X/Y}(\log\Delta)\]
We define the Hodge filtration of $V $ as
\[\{F^p\}=\{\text{Im}(R_\text{dR}^kf_* F_\text{tru}^p\Omega^\bullet_{X/Y}(\log\Delta)\to R_\text{dR}^k f_*\Omega^{\bullet}_{X/Y}(\log S))\}\subset V.\]
The connecting maps on $R^i_\text{dR}f_*$ of
\\[.2cm]
{\bf (2.1.5.)}
\[ 0\to f^*\Omega^1_Y(\log S)\otimes F_\text{tru}^{p-1}\Omega^{\bullet}_{X/Y}(\log\Delta)[-1]\to
F_\text{tru}^pK^0/K^2\Omega^\bullet_X(\log\Delta)\]
\[\to F_\text{tru}^p\Omega^{\bullet}_{X/Y}(\log\Delta)\to 0\] is compatible
with the Gauss-Manin connection under the natural map induced by the inclusion of short exact sequences ${\bf (2.1.5)\subset (2.1.4.)}$, and
we find\\[.2cm]
{\bf Griffiths' transversality}
$$\nabla F^p\subset\Omega^1_Y(\log S)\otimes F^{p-1},$$
which gives rise to a filtered logarithmic de Rham bundle
$$(V,\nabla, F^\bullet),$$
and with the associated graded Higgs bundle
$$\text{Gr}_{F^\bullet}(V,\nabla)=\left(\bigoplus_{p=0}^k{F^p\over F^{p+1}},\quad\bigoplus_{p=0}^k{F^p\over F^{p+1}}\xrightarrow{\bar\nabla}\Omega^1_Y(\log S)\otimes{F^{p-1}\over F^{p}}\right).$$
As over complex number field $E_1$-degeneration of Hodge to de Rham holds true. It implies that the map induced by the inclusion of the complexes ${\bf (2.1.5.)\subset (2.1.4.)}$ is an isomorphism
\[R_\text{dR}^kf_* F_\text{tru}^p\Omega^\bullet_{X/Y}(\log\Delta)\simeq F^p.\]
Consequently it induces an isomorphism
\[(\bigoplus_{p+q=k}E^{p,q},\bigoplus_{p+q=k}\theta^{p,q})\simeq
\text{Gr}_{F^\bullet}(V,\nabla).\]
By Deligne the filtered logarithmic de Rham bundle $(V,\nabla, F^\bullet)$ is the quasi-canonical extension of the
variation of polarized variation of Hodge structure of $k$-th Betti cohomology
of the smooth family $$f:V\to U.$$\\[.4cm]
\subsection{ Positive (Negative) Sheaves in Hodge Theory}  There are several different notions of positivities of vector bundles and sheaves.\\[.2cm]
{\bf Definition 2.2.1.}(Griffiths) A holomorphic vector bundle $E$ over a complex manifold $Y$ called semi-positive if $E$ carries a Hermitiam metric whose curvature form ((1,1)-type and Hermitian)
$$\Theta_E(e,\xi)\geq 0$$
for any section $e\in E$ and any holomorphic tangent vectot $\xi$.
\\[.2cm]
In algebraic Setting:\\[.2cm]
{\bf Definition 2.2.2.}(Viehweg) Let $\mathcal F$ be a torsion free sheaf on a quasi-projective normal
variety $Y$, $\mathcal{H}$ be an ample invertible sheaf and $U\subset Y$ be an open
subvariety.
\begin{itemize}
\item
[1] $\mathcal F$ is globally generated over $U$ if the natural morphism
$$ H^0(Y,\mathcal F)\otimes\mathcal O_Y\to\mathcal F$$ is surjective over $U.$
\item
[2] $\mathcal F$ is weakly positive over $U$ if the restriction of $\mathcal F$ to $U$ is locally
free and if for each $\alpha>0$ there exists some $\beta>0$ such that
$$ S^{\alpha\beta}(\mathcal F)\otimes\mathcal H^\beta$$ is globally generated over
$U$.
\item
[3] $\mathcal F$ is ample w.r.t. $U$ if there exist some $\eta>0$ and a morphsim
$\bigoplus\mathcal H\to S^\eta(\mathcal F)$
is surjective over $U.$
\end{itemize}
In short we just call $\mathcal F$  weakly positive in 2), or generically ample in 3) without refereeing an open subvariety.
\\[.2cm]
{\bf Property 2.2.3. (of weakly positive sheaves)}
\begin{itemize}
\item [a] The notion weakly positivity does not depend on the choice of an ample invertible sheaf and
does not depend on base change unter generically finite morphism.
\item [b]If $\mathcal F$ is weakly positive and if $\mathcal F\to\mathcal G$ is surjective on an open set of $Y$ then $\mathcal G$ is weakly positive.
\item [c] A tensor product of weakly positive sheaf with an ample sheaf is ample.
\item [d] Assume $Y$ is projective, a numerically effective vector bundle $\mathcal F$ ( i.e. any quotient bundle of $\mathcal F|_C$ restricted to any projective curve
$C\subset Y$ has non-negative degree) is weakly positive.
\end{itemize}.\\[.2cm]
{\bf Examples of semi-positive Hodge bundles}.\\[.2cm]
 Let $$ (E,\theta)=(\bigoplus E^{p,q},\bigoplus\theta^{p,q})$$ be a system of Hodge bundles arising from the canonical extension of a polarized variation of Hodge structure over $(Y, S).$\\[.2cm]
 Griffiths: (for $S=\emptyset$)
\begin{itemize}
\item $E^{k,0}$ is semi-positive. Hence, it
is weakly positive on  $Y.$
\item The Griffiths' augmented line bundle
$$\Lambda:=(\det E^{k,0})^k\otimes(\det E^{k-1,1})^{k-1}\otimes\cdots\otimes (\det E^{0, k-1})$$
is semi-positive, and it is strictly positive on an open set if the Torelli map is generically injective.
\end{itemize}
{Kawamata:} The above results hold true over a log base.\\[.2cm]
Griffiths' curvature formular for Hodge metric and, in general, Simpson's semi-stability for Higgs bundles carrying Yang-Mills-
Higgs metric allow us to study the negativity of the the kernel of the Higgs field:
$$K^{p,q}=\text{Ker}( E^{p,g}\xrightarrow{\theta^{p,q}}E^{p-1,q+1}\otimes\Omega^1_Y(\log S)).$$
\\[.2cm]
{\bf Proposition 2.2.4.} (Zuo, Brunebarbe)  {\sl The dual ${K^{p,q}}^\vee$ is weakly positive.}
\\[.2cm]
{\bf Sketch the proof.}
$K^{p,q}\subset E^{p,q}$ is torsion free subsheaf, and after a blowing up $Y'\to Y$ on the singular locus of $K^{p,q}$ we may assume it is locally free.
We show
${K^{p,q}}^\vee $ is nef by\\[.2cm]
$\bullet$ Griffiths' curvature formular for the Hodge metric
on $U$.\\[.2cm]
$\bullet$ The works by Schmid, Cattani-Kaplan-Schmid,
Kashiwara and Steenbrink on
variation of mixed Hodge structure along the boundary divisor $ S.$
The Hodge metric $h$ has singularity along $ S$ of logarithmic growth. Hence
the curvature current $\Theta_{(K^{p,q}, h)}$ represents the Chern classe of $K^{p,q}$ on $Y$.\\[.2cm]
Griffiths's curvatura formular reads
$$\Theta_{(E^{p,q}, h)}=-{\theta^{p,q}}^*\wedge\theta^{p,q}-\theta^{p-1, q+1}\wedge{\theta^{p-1,q+1}}^*$$
$$=\text{Positive}+\text{Negative}.$$
Restricting tha to $K^{p,q}$ we find a consulation of the positive term $$({\theta^{p,q}}^*\wedge\theta^{p,q})|_{K^{p,q}}=0,$$
hence
$$\Theta_{( K^{p,q},h)}=-(\theta^{p-1, q+1}\wedge{\theta^{p-1,q+1}}^*)|_{K^{p,q}}+\sigma^*\wedge\sigma\leq 0,$$
where $\sigma$ is the second fundamental form of $K^{p,q}\subset E^{p,q}.$\\[.2cm]
It shows:\\[.2cm]
{\bf Special Case.} For a projective testing cuve $C\subset Y$ with $C\not\subset S$ one has
${K^{p,q}}^\vee|_C$ is semi-positive.\\[.2cm]
In general, for a curve $C\subset S\subset Y$ we find the minimal intersection
of the components of $\bar S$
$$ C\subset S_I\quad\text{and}\quad C\not\subset S_J\quad\forall J\supset I.$$
The multi residue map $\text{Res}_I\theta$ on $E|_{ S_I}$
along the components $S_i,\, i\in I$ defines the weight filtration on $E|_{ S_I}$, and it restricts to
a filtration on $K^{p,q}|_{ S_I}$
$$W_{0,\cdots 0}(K^{p,q}|_{ S_I})\subset\cdots\subset W_{2k,\cdots, 2k}(K^{p,q}|_{ S_I})=K^{p,q}|_{S_I},$$
such that quotients lie in the kernels of the Higgs fields on the graded Higgs bundles of
the polarized VHS arising from the grading of the weight filtration on the variation of mixed Hodge structures along $ S.$ The proof follows from the special case.\\[.4cm]
We state some immediate conequences from the weakly positivity of ${K^{p,q}}^\vee.$\\[.2cm]
{\bf Proposition 2.2.5.} {\sl If $U$ carries a generically injective period mapping $$\phi: U\to\mathcal D. $$ Then
the logarithmic canonical line bundle $\omega_{ Y}( S)$ is big, i.e. ample w.r.t. an open set of $ Y.$}\\[.2cm]
{\bf Sketch the proof}:\\[.2cm]
{\bf 1.} $\Omega^1_{ Y}(\log S)$ is weakly positive.\\[.2cm]
The Higgs field $\theta: E\to E\otimes\Omega^1_{ Y}(\log S)$ induces a natural map of Higgs bundles
$$d\phi: (T_{ Y}(-\log S), 0)\to\mathcal End(E,\theta) $$
which coincides with the derivative of the period mapping.\\[.2cm]
As the image $d\phi T_{ Y}(-\log S)=:T'$ lies in the kernel the Higgs field on $\mathcal{E}nd(E)$ we have 
$${T'}^\vee\subset\Omega^1_{Y}(\log S)$$ is weakly positive.
And as ${T'}^\vee$ coincides with $\Omega^1_{ Y}( S)$ on an open sub set, $\Omega^1_{ Y}(\log S)$
is weakly positive.\\[.2cm]
{\bf 2.} $\omega_{\bar Y}(\bar S)$ is big.\\[.2cm]
Consider the Griffiths's augmented line bundle $\Lambda $ on $Y$ which is big as the period mapping is injective at a point. One takes a suitable tensor product the original Higgs bundle
$$(E_T,\theta_T)= (E,\theta)^{\otimes m_0}\otimes\cdots\otimes (E,\theta)^{\otimes m_{k-1}}$$
such that $\Lambda$ is contained a Hodge bundle $ E_T^{r_0,s_0}$. By runing the maximal non-zero iteration of the Higgs map $\theta_T$ with the initial line bundle $\Lambda\subset E_T^{r_0, s_0}$ we obtain the Griffiths-Yukawa coupling
$$S^\ell T_Y(-\log S)\otimes\Lambda\xrightarrow{\theta_T^\ell\not=0} E_T^{r_0-\ell, s_0+\ell},$$
such that the image of $\theta_T^\ell$ lies in $K=\text{Ker}(\theta_T^{r_0-\ell, s_0+\ell}).$
The dual gives rise to a big sub sheaf in a symmetric power
of log differential forms and with a weakly positive quotient:
$$\Lambda\otimes K^\vee\xrightarrow{{\theta_T^\ell}^\vee\not=0}S^\ell\Omega^1_Y(\log S)\to Q\to 0.$$
By taking determinants for the above exact sequence we obtain
$$\omega(S)\supset\Lambda^\alpha\otimes P^\beta,$$
where $P$ is a weakly poisitive line bundle and $\alpha,\beta >0$.\\[.4cm]
\subsection{ Graded Higgs Bundle Arising From Kodaira-Spencer Deformation Theory}
Given a log smooth family $$f: (X,\Delta)\to (Y, S).$$   Replaceing the injectivity of Torelli mapps attached to $f$ by the non-isotriviality of the family
 Viehweg-Zuo showed that the non-linear variation of complex structure on the fibres still produces certain negativities on the log base. \\[.2cm] We start with the original (logarithmic)  Kodaira-Spencer map 
$$ T_Y(-\log S)\xrightarrow{\tau^{n,0}}R^1f_*T_{X/Y}(-\log\Delta),$$
and define then the extended (logarithmic) Kodaira-Spencer map   $\tau^{p,q}$ by 
\begin{equation*}
\xymatrix@R=2cm{
T_Y(-\log S)\otimes R^qf_*T^q_{X/Y}(-\log\Delta)\ar[d]^-{\tau^{n,0}\otimes\mathrm{Id}}\ar[r]^-{\tau^{p,q}}& R^{q+1}f_*T^{q+1}_{X/Y}(-\log\Delta)\\
R^1f_*T_{X/Y}(-\log\Delta)\otimes R^qf_*T^q_{X/Y}(-\log\Delta).\ar[ur]^-{\cup}
}
\end{equation*}
Putting all individual sheaves
$$F^{p,q}=R^qf_*T^q_{X/Y}(-\log\Delta)/\text{torsion},\quad p+q=n$$
$$\left(= R^qf_*(\Omega^p_{X/Y}(\log\Delta)\otimes
(\Omega^n_{X/Y}(\log\Delta))^{-1})/\text{torsion}\right),$$
and the maps $\tau^{p,q}$  together
we obtain the so-called {\sl Deformation Higgs bundle (sheaf)} 
$$(F,\tau):=(\bigoplus_{p+q=n}F^{p,q},\bigoplus_{p+q=n}\tau^{p,q})$$
attached to $$f:X\to Y.$$
One checks that $\tau=\oplus\tau^{p,q}$ satisfies the integrability condition $\tau\wedge\tau=0$ using the associativity and anti-commutivity of the cup product on Dolbeault cohomology of the tangent sheaves along fibres.
The Kodaira-Spencer map on the deformation Higgs bundle   is tautologically non-trivial for a non-isotrivial family. The crucial property is that
{\sl the kernel of the Kodaira-Spencer map is strictly negative,}
 which follows from a comparison theorem
between the deformation Higgs bundle and a system of Hodge bundles and Kawamata-Viehweg's theorem on the positivity of direct image sheaves.\\[.2cm]
As a motivation we consider a semi-stable family $$f: X\to Y$$ of abelian varieties, or Calabi-Yau $n$-folds with the maximal Var($f$), and find the relation between the deformation Higgs bundle and the Higgs bundle arisng from VHS on the middle cohomology of fibres in the form
$$(F,\tau)=(E,\theta)\otimes{E^{n,0}}^{-1},$$
where the Hodge line bundle $E^{n,0}$ is ample on $Y$.  It is clear that   from the above comparison map the kernel of $\tau$ is strictly negative.\\[.2cm]
In general, for a family $$f:X\to Y$$ of varieties with good minimal models, or of general type the theorem due to Kawamata, Viehweg and Koll\'ar on the positivity of the direct image sheaves 
of the relative pluri canonical sheaf $f_*\omega_{X/Y}^\nu$ allows us in the next section to construct a non-trivial Higgs map
from $(F,\tau)$ to a system of Hodge bundles $(E,\theta)$ arising from an intermediate geometric situation
$$\rho: (F,\tau)\to (E,\theta)\otimes A^{-1}$$
 after taking a higher power of the self-fibre product of the orginal family.  
\\[.4cm]
\section{\bf Comparison Between Deformation Higgs Bundle and System of Hodge Bundle}
\subsection{Positivity of Direct Image Sheaves}
Let $f:X\to Y$ be a log smooth family, whose reduced singular fibres $\Delta$  is normol crossing over $S$ and
$$f: V=X\setminus\Delta\to Y\setminus S=:U$$
 the smooth part of the family.
The following positvity of $f_*\omega_{X/Y}^\nu$ plays the fundamental role in study Iitaka conjecture:\\[.2cm]
{\bf Theorem 3.1.1.}
(Kawamata and Viehweg) {\sl Assume that $f: X\to Y$ has the maximal variation, and $\omega _{V/U}$ is semi ample. Then $f_*\omega_{X/Y}^\nu $ is ample w.r.t. open subsets for all $\nu>1$ with $f_*\omega_{X/Y}^\nu\not=0.$}
.\\[.2cm]
For constructing the comparison
map we need actually a slightly stronger form of the positivity for the reletive log top differential forms
$$\mathcal{L}:=\Omega^n_{X/Y}(\log\Delta)\subset\omega_{X/Y}.$$
\\[.2cm]
{\bf Proposition 3.1.2.}(Viehweg-Zuo)
\\[.2cm]{\sl
{\bf 1.}
$f_*(\mathcal L^\nu)$ is is ample on an open subset $U'\subset U$ for $\nu\gg 0$.\\[.2cm]
{\bf 2.} Fixing an ample line bundle $A$ on $ Y$ and replacing the original family $f: X\to Y$ (a birational modification) by a self fibre product of a higher power $f^{(r)}: X^{(r)}\to Y$,
or by a Kawamata base change $\phi: Y'\to Y$ we may assume $\mathcal L^\nu\otimes f^*(A^{-\nu})$ is globally generated over $f^{-1}(U').$}
\\[.2cm]
\subsection{Comparion Map}.\\[.2cm]
{\bf Hodge Theory on Cyclic Covers.} One of the motivation of constructing cyclic cover is trying to give a more Hodge theoretical proof for Kodaira-Akizuki-Nakano vanishing theorem.\cite{EV}\\[.2cm]
Let $X$ be a projective manifold, $\mathcal L$ a line bundle on $X$ and $s\in H^0(X,\mathcal L^\nu)$ with the zero divisor $D\subset X$.  For the simplicity, we assume $D$ is smooth.  One takes the $\nu$-th cyclic cover 
$\gamma: Z\to X$
with 
$$\gamma_*\Omega^p_Z(\log \gamma^*D)
=\bigoplus_{i=0}^{\nu-1}\Omega^p_X(\log D)\otimes \mathcal L^{-i}.$$
Deligne has shown
$$ H^k(Z\setminus\gamma^*D,\mathbb C)=\bigoplus_{p+q=k} H^q(Z,\Omega^p_Z(\log \gamma^*D)).$$
Assume $D$ is ample, then $X\setminus D$ is affine, the same holds true for $Z\setminus \gamma^*D$
and hence $$H^k(Z\setminus \gamma^*D, \mathbb C)=0,\quad  \forall\, k>\dim X=n.$$\\
By the Hodge decomposition 
$$0=H^q(Z,\Omega^p_Z(\log\gamma^* D)=
\bigoplus_{i=0}^{\nu-1}H^g(X, \Omega^p_X(\log D)\otimes\mathcal L^{-i}).$$
for any $p+q>n.$  In particular,
$$ H^q(X,\Omega^p_X(\log D)\otimes\mathcal L^{-1})=0,\forall\, p+q>n.$$
Using the residue map and induction on $\dim X$ we
show Kodaira-Akizuki-Nakano vanishing theorem
$$ a) \quad H^q(X,\Omega^p_X\otimes\mathcal L)=0,\quad \forall\, p+q>n$$
$$ b) \quad H^q(X,\Omega^p_X\otimes\mathcal L^{-1})=0,\quad \forall\, p+q<n.$$
\\[.4cm]
We remark, in contrast, the middle dimensional cohomology $H^q(X,\Omega^p_X\otimes \mathcal L^{-1}),\quad p+q=n$
will not be zero in general and play the crucial role in the construction of the comparsion map.\\[.2cm]
{\bf Construction 3.2.1.}
For $\mathcal L=\Omega^n_{X/Y}(\log\Delta),$  by 2 in Prop.3.1.2. after taking a power of the self fibre product we may assume $\mathcal L^\nu\otimes f^*(A^{-\nu})$ is globally generated over $f^{-1}(U')$  for some open set in $U.$ We 
choose a generic non-zero section $s$ of $\mathcal L^\nu\otimes f^*A^{-\nu}.$ Then the $\nu$-th cyclic cover defined by $s$ induces a family $$g: Z\to X\xrightarrow{f}Y$$ with the singular fibres $\Pi$ over $S+T$, where $T$ is the degeneration locus of the \emph{new} singular fibers arising from the cyclic cover. By blowing up and leaving out a codimension-two sub scheme of $Y$ we may assume
$$g: (Z,\Pi)\to (Y, S+T) $$ is log smooth.\\[.2cm]
Take $(E,\theta)$ to be the graded Higgs bundle of Deligne's quasi-canonical extension of VHS on the middle cohomology $ R^n_\text{B}g_*\mathbb Z_{Z\setminus\Pi}$
on $Y,$
then there exists a Higgs map $$\rho: (F,\tau)\to (E,\theta)\otimes A^{-1},$$
that is
$$
\xymatrixcolsep{5pc}\xymatrix{
F^{p,q}\ar[r]^-{\tau^{p,q}}\ar[d]^{\rho^{p,q}}& F^{p-1,q+1}\otimes\Omega^1_{{Y}}(\log {S})\ar[d]^{\rho^{p-1,q+1}\otimes\iota}\\
A^{-1}\otimes E^{p,q}\ar[r]^-{\mathrm{id}\otimes\theta^{p,q}}& A^{-1}\otimes E^{p-1,q+1}\otimes\Omega^1_{\bar{Y}}(\log ({S}+{T})).
}
$$
where  $$\iota: \Omega^1_{ Y}(\log S)\hookrightarrow \Omega^1_{ Y}(\log( S+ T))$$
 is the natural inclusion.\\[.2cm]
Note that the comparison map $\rho^{p,q}$ is defined only on $Y$ \emph{a priori}, that is, a morphism between $F^{p,q}=\bar F^{p,q}|_Y$ and $E^{p,q}=\bar E^{p,q}|_Y$. Since $\bar F^{p,q}$ is reflexive and $\mathrm{codim}(\bar{Y}\setminus Y)\geq 2$, the comparison map $\rho^{p,q}$ extends to $\bar Y$.\\
We emphasize the crucial point in the comparison map:
 although the Higgs field $\theta$ on $E$ has singularity along $\bar S+\bar T$, its restriction to $\rho(F)$ has only singularity on the original degeneration locus $\bar S.$\\[.2cm]
We give more detailed discussions on the construction 3.2.1.  below:\\[.2cm]
{\bf Step 1. Blowing Up: $X'\to X$.}\\[.2cm]
Choose a generic non-zero section $s$ of $\mathcal L^\nu\otimes f^*(A^{-\nu})$
such that the zero divisor $D:=\mathrm{div}(s)$ intersects the generic fibres of $f$ smoothly.
Let $T\subset Y$ denote the closure of the discriminant of the map
$$f: D\cap V\to U,$$ which is the locus where the intersection
$$D_y:= D\cap f^{-1}(y)$$ becomes singular. \\[.2cm] Let $\Sigma:=f^{-1}(T),$ we 
take a blowing up $$\delta: X'\to X$$ with centers in $D+\Delta+\Sigma.$ 
Leaving out some codimension-2 subscheme,
we may assume that $\delta^*(D+\Delta+\Sigma)$ is normal crossing and the composition map
$$f': X'\xrightarrow{\delta}X\xrightarrow{f}Y$$
is log smooth morphism between the log pairs
$$ f' :\left(X',\delta^*(D+\Delta+\Sigma)\right)\to (Y, (S+T)).$$
\\[.4cm]
{\bf Step 2. Cyclic Cover Defined by $s$.}\\[.2cm]We write $\mathcal M:=\delta^*(\mathcal L\otimes f^*A^{-1})$ and $D':=\delta^*(D)$, then
$\mathcal M^\nu=O_{X'}(D')$. One takes the $\nu$-th cyclic cover for the section $\delta^*s\in H^0(X',\mathcal M^\nu)$ and the normalization
$$\gamma': Z'\xrightarrow{\text{normalization}}X'(\sqrt[\nu]{\delta^*s})\xrightarrow{\gamma}X'.$$
$Z'$ could be singular.
By taking a resolution of singularity of $ Z'$, and a blowing up at the centers in the fibres
over $Y$ we obtain a non-singular variety $Z$ and a birational morphism
$\eta: Z\to Z'$. Againg leaving out a codimension-2 subscheme
in $Y$ we may assume the induced map
$$g:Z\xrightarrow{\eta}Z'\xrightarrow{\gamma'}X'\xrightarrow{f'}Y$$
is log smooth for the pairs
$$g: (Z , g^{-1}(S+T))\to (Y, (S+T)).$$
Let $$g: Z_0\to U$$ be the smooth part of the map $g.$\\[.4cm]
{\bf Step 3. Differential Forms On the Cyclic Cover}.\\[.2cm]
Recall that the local system $\mathbb V=R^n_\text{B}g_*\mathbb C_{Z_0}$ gives rise to the filtered logarithmic de Rham bundle
$$\nabla: V\to V\otimes\Omega^1_Y(\log (S+T))$$
as the quasi canonical extension
of $\mathbb V\otimes\mathcal O_{Y\setminus(S+T)},$ and with 
 the induced system of Hodge bundles
$$\text{Gr}_F(V,\nabla)=(\bigoplus_{p+q=n}R^qg_*\Omega^p_{Z/Y}(\log\Pi),\bigoplus_{p+q=n} \theta^{p,q})          =:(\bigoplus_{p+q=n}E^{p,q},\bigoplus_{p+q=n}\theta^{p,q})=(E,\theta).$$
The Higgs map
$$\theta^{p,q}: E^{p,q}\to E^{p-1,q+1}\otimes\Omega^1_Y(\log (S+T))$$
is
the edge map of $R^\bullet g_*$ of the exact sequence\\[.2cm] {\bf (3.2.1)} $$ 0\to g^*\Omega^1_Y(\log(S+T))\otimes\Omega_{Z/Y}^{p-1}(\log\Pi)\to
\text{gr}\Omega^p_Z(\log\Pi)\to
\Omega^p_{Z/Y}(\log\Pi)\to 0.$$
Via the blowing up $$\delta: X'\to X$$
we consider now the pulled back of the deformation Higgs bundle $(F,\tau)$ on $Y$
$$\delta^*(F,\tau)=(\bigoplus\delta^*F^{p,q},\bigoplus\delta^*\tau^{p,q})=(\bigoplus F'^{p,q},\bigoplus\tau'^{p,q}) $$
with
$$F'^{p,q}=R^qf'_*(\delta^*\Omega^p_{X/Y}(\log\Delta)\otimes\delta^*\mathcal L^{-1})/\text{torsion}.$$
And note that the Kodaira-Spencer map
$$\tau'^{p,q}: F'^{p,q}\to F'^{p-1,q+1}\otimes\Omega^1_Y(\log S)$$
as the edge map of $R^\bullet f'_*$ of the exact seuqence\\[.2cm]
{\bf (3.2.2)}
$$0\to f'^*\Omega^1_Y(\log S)\otimes\delta^*\Omega^{p-1}_{X/Y}(\log\Delta)\otimes\mathcal L'^{-1}\to
\delta^*\text{gr}\Omega^p_X(\log\Delta)\otimes\mathcal L'^{-1}\to\delta^*\Omega^p_{X/Y}(\log\Delta)\otimes\mathcal L'^{-1}\to 0.$$\\[.4cm]
{\bf Step 4 Comparison Between two Higgs Bundles}\\[.2cm]
Let $\bullet$ stand either for $\text{Spec}(\mathbb C)$ or for $Y$. Then the Galois group $\mathbb Z/\nu\mathbb Z$ of
$$\psi: Z\to X'$$ acts on
$\psi_*\Omega^p_{Z/\bullet}(\log\Pi)$
with the eigenspace decomposition
$$\psi_*\Omega^p_{Z/\bullet}(\log\Pi)
=\Omega^p_{X'/\bullet}(\log(\Delta'+\Sigma'))\oplus
\bigoplus_{i=1}^{\nu-1}\Omega^p_{X'/\bullet}(\log(\Delta'+\Sigma'+D')\otimes\mathcal L'^{-1}\otimes f'^*A,$$
which indues a natural inclusion\\[.2cm]
{\bf (3.2.3)}
\begin{equation*}
\xymatrix@R=2cm{
\delta^*\Omega_{X/\bullet}^p(\log\Delta)\otimes\mathcal L'^{-1}
\ar@{^(->}[r]
\ar@{^(->}[d]
&
\psi_*\Omega_{Z/\bullet}^p(\log\Pi)\otimes f'^*(A^{-1})\\
\delta^*\Omega_{X/\bullet}^p(\log\Delta+\Sigma)\otimes\mathcal L'^{-1}
\ar@{^(->}[r]
&
\Omega_{X'/\bullet}^p(\log\Delta'+\Sigma'+D')\otimes\mathcal L'^{-1}\otimes f'^*(A)\otimes f'^*(A^{-1})
\ar@{^(->}[u]\\
}
\end{equation*}
and via $\psi: Z\to X'$ the inclusion {\bf (3.2.3)}   together with the inclusion
$$\Omega^1_Y(\log S)\hookrightarrow\Omega^1_Y(\log (S+T))$$
induce an inclusion of the exaxt sequences
$$\psi^*{\bf (3.2.2)}\subset{\bf (3.2.1)}\otimes g^*A^{-1}.$$
Finally the direct image $$g_*\left(\psi^*{\bf (3.2.2)}\subset{\bf (3.2.1)}\otimes g^*A^{-1}\right)$$
it yields a map
$$\rho^{p,q}: F^{p,q}\to E^{p,q}\otimes A^{-1}$$
commuting with the Higgs fields as the edge maps.\\[.4cm]
\subsection{ Maximal Non-Zero Iteration of Kodaira-Spencer Map.}  We need to show the non-triviality of the comparison map.\\[.2cm]
{\bf Proposition 3.3.1.} (Torelli type injectivity){\sl
\begin{itemize}
\item[1] $\rho^{n,0}: O_{\bar Y}=F^{n,0}\xrightarrow{s}E^{n,0}\otimes A^{-1}$ is injective.
\item [2.] If the canonical line bundle of the fibres is ample, then every $\rho^{p,q}$ is injective at all points of $U\setminus T$.
\item [3] If the canonical line bundle of the fibres is big, then $\rho^{n-1,1}$ is pointwisely injective in $U\setminus T.$
\item [4] If the canonical line bundle on the fibres is semi-ample, then the following composition map
$$ T_{\bar Y}(-\log\bar S)\otimes F^{n,0}\xrightarrow{\rho^{n,0}}T_{\bar Y}(-\log\bar S)\otimes\rho^{n,0}(F^{n,0})\xrightarrow{\theta^{n,0}}E^{n-1,1}\otimes A^{-1}$$
is pointwisely injective in a Zariski open subset of $U\setminus T.$
\end{itemize}}.\\[.2cm]
Sketch the proof of Prop. 3.3.1.:
\begin{itemize}
\item[2] follows from Akizuki-Kodaira-Nakano vanishing theorem.
\item[3] follows from Bogomolov vanishing theorem.
\item[4] Viehweg-Zuo \cite{VZ-1} showed a weaker form of (4):
the map 
$$T_Y(-\log S)\otimes F^{n,0}\xrightarrow{\theta^{n,0}\circ\rho^{n,0}}E^{n-1,1}\otimes A^{-1}$$
is non-zero simply by using the global form of the semi-negativity of the kernel of the Higgs map.\\
Deng \cite{Deng}  obtained the full statement for (4)  basically using
the Griffiths curvature formular and the the positivity of the curvature of $A$ pointwisely.
\end{itemize}.\\[.2cm]
For families of hypersurfaces Carlson, Green anf Griffiths studied iterated cup product with Kodaira-Spencer map on the systems of Hodge bundles from the middle cohomologies.\\
For the case of Calabi-Yau manifolds the importance of this coupling was brought up by physicists. So, 
we shall call it as the Griffiths-Yukawa coupling.\\
Viehweg-Zuo adopted this notion for deformation of Higgs bundle for producing a big subsheaf in a symmetric power of $\Omega^1_{\bar Y}(\log\bar S)$  as follows.\\[.2cm]
{\bf Proposition 3.3.2.} (Viehweg-Z) 
{\sl Some symmetric power of the log cotangent sheaf on the base contains an ample sub sheaf
	$$  \mathcal A\hookrightarrow S^\ell\Omega^1_{\bar Y}(\log\bar S).$$}\\[.2cm]
{\bf Sketch the proof:}
We consider the sub Higgs sheaf
$$(G,\tau_G):=\rho ((F,\tau))\subset (E,\theta)\otimes A^{-1},$$
and note the crucial property:\\[.2cm]
{\bf 1}
$G^{n,0}\simeq\mathcal O_{\bar Y}$ and
$\tau_G$ has logarithmic pole only along the original boundary divisor $\bar S.$ \\[.2cm]
{\bf 2}
The Torelli type injectivity ( Prop.3.3.1.)
shows that in any case we have non-zero Higgs map
$$\tau_G^{n,0}: T_{\bar Y}(-\log\bar S)\to G^{n-1,1}.$$
 Which allows us to run the maximal non-zero iteration of the Higgs map on $(G,\tau_G) $ in the form:
 $$S^\ell T_{\bar Y}(-\log\bar S)\xrightarrow{\tau_G^\ell\not=0} G^{n-\ell,\ell}\subset E^{n-\ell,\ell}\otimes A^{-1}$$
 and such that $\tau_G^{\ell+1}=0$. Hence, the map factors through
 $$ S^\ell T_{\bar Y}(-\log\bar S)\xrightarrow{\tau_G^\ell\not=0}\text{ker}(\tau_G^{n-\ell,\ell})\subset\text{ker}(\theta^{n-\ell,\ell})\otimes A^{-1}.$$
\\[.2cm]
As      $\text{ker}(\theta_G^{n-\ell,\ell})\leq 0$ the dual of the Griffiths-Yukawa coupling gives rise to a generically ample subsheaf
 $$\mathcal A:=\text{im}\left(A\otimes\text{ker}(\theta^{n-\ell,\ell})^\vee\xrightarrow{{\tau_G^\ell}^\vee\not=0}
 S^\ell\Omega^1_{\bar Y}(\log\bar S)\right).$$\\[.2cm]
\subsection{ Second Construction of The Comparison Map via Kawamata's Cover.}  
In the first construction of the comparion map in Section 3.2.  we use the higher power of self fibre product of the original family and we have no control on the power $\ell$ such that $S^\ell\Omega^1_{\bar Y}(\log\bar S)$ containing a
big sub sheaf. Here we shall briefly discuss on the second contruction using Kawamata's cover.\\[.2cm]
$\bullet$ Fixing an ample line bundle such that $A\subset \mathcal L^\nu$
for an $\nu>>0$.
By Kawamata there exists a Kummer cover
$\psi: Y'\to Y$ so that $\psi^*A=A^{'\nu}$ for an ample line bundle $A'$ on $Y'.$\\[.2cm]
$\bullet$ Let $f': X'\to Y'$ denote the fibre prodcut of the base change, 
then $\mathcal L'=f'_*\Omega^n_{X'/Y'}(\log\Delta')$ and
$$\mathcal L^{'\nu}=\mathcal\psi^*\mathcal L^\nu\supset A^{'\nu}.$$

By taking the $\nu$-th cyclic cover $\pi': Z'\to X'$ of the canonical section $s$ of $(\Omega^n_{X'/Y'}(\log\Delta')\otimes f^{'*}A'^{-1})^\nu$.
\\[.2cm]
$\bullet$ Applying the construction 3.2. for a section of $\mathcal L^{'\nu}\otimes A^{'-\nu}$ on $Y'$ we obtain a Higgs map
$$\rho:\psi^*(F,\tau)\to (E,\theta)\otimes A^{'-1},$$
where $(E,\theta)$ is the system of Higgs bundles of the quasi-canonical extension of a PVHS of weight-$n$ over $Y'\setminus (\psi^*(S)+T')$, where $n$ is the dimension of the fibres of the original family $f:X\to Y.$
We obtain a big subsheaf
$\mathcal A\hookrightarrow S^\ell\psi^*\Omega^1_Y(\log S) $
with $\ell\leq n.$\\[.4cm]
\subsection{ A Remark to Families of $n$-Folds of General Type or of Good Minimal Model.}.\\[.2cm]
For a family of $n$-folds either of good minimal model or of general type
$$ f:\bar X\to \bar Y$$ Then Kawamata for good minimal model \cite{K} Kollar for general type \cite{Kol}
showed that $f_*\omega_{\bar X/\bar Y}^\nu$ is big for $\nu>>0$.\\[.2cm]
The main difference between the case of good minimal model and the case of semi-ample is that the linear system of $\omega_{\bar X/\bar Y}^\nu$ in the first case could be not globally generated over $f^{-1}(U_0)$ for any open subset of $U$, while in the second case it is globally generated over $f^{-1}(U_0)$ for some open subset of $U.$\\[.2cm]
Popa-Schnell applied the theory of Hodge module to get a comparison
similar to what stated in Section 3.2. It has the advantage that one does not care about the complication of the singularity appearing in the construction (of codimension-2)  too much.\\
Below we propose an approach along the original construction in Section 3.2., which already appeared in \cite{VZ-1} over 1-dimensional base.\\
After taking a higher power of self-fibre product of $f$ one finds a non-zero section $s$ of $({\Omega^n_{X/ Y}(\log\Delta)\otimes A^{-1}})^\nu$ with the zero divisor $D$ of $s.$ Note that the intersection of $D$ with the generic fibre could be singular.\\[.2cm]
$\bullet$ Replacing the family $$f: X\to Y$$ by a blowing up
$$f': X'\xrightarrow{\delta}X\to Y$$
such that $\delta^*(D)=:D'$
is a normal crossing divisor.\\[.2cm]
$\bullet$ Let $T'\subset Y$ denote the closure of the discriminant of the map
$$f: D'\cap V'\to U, $$ that is the locus where the normal crossing divisor
$D'$ meets $f^{'-1}(y),\, y\in U$ not transversally.\\[.2cm] Let
$\Sigma'=f^{'-1}(T')$, 
take a further blowing up
$$f'': X''\xrightarrow{\delta'}X'\xrightarrow{f'}Y$$
and leave out some codim-2 sub scheme in $Y$
such that $f''$
is a log smooth morphism between the log pairs
$$ f'' :\left(X'',\delta^{'*}(D'+\Delta'+\Sigma')\right)\to (Y, (S+T')).$$
$\bullet$ Take the $\nu$-th cyclic cover
$$X''(\sqrt[\nu]{\delta^{'*}\circ\delta^*(s)})\xrightarrow{\gamma}X'',$$
and we procedure further as in Section 3.2.
\\[.4cm]
\section{ \bf Algberaic Hyperbolicity and Viehweg Hyperbolicity of Moduli Stacks}
\subsection{ Algebraic Hyperbolicity on Moduli Stacks.}
Given a family $f: V\to U$ of $n$-folds with semi-ample canonical line bundle and such that the classifying map from $U$ to the moduli space is quasi-finite.
We take a smooth projective compactification $\bar Y\supset U$ such that $\bar S=\bar Y\setminus U$ is a normal crossing divisor.\\[.2cm]
As an application of the existence of a big sub sheaf in $S^\ell\Omega^1_{\bar Y}(\log \bar S)$ we study the the algebraic hyperbolicity $(\bar Y,\bar S)$ and give a uniform proof for the result on algebraic hyperbolicity of $(\bar Y,\bar S)$ due to Parshin-Arakelov, Miglinorini, Zhang, Kov\'acs, Bedulev-Viehweg, Oguiso-Viehweg and Viehweg-Zuo.\\[.2cm]
{\bf Proposition 3.1} {\sl Assume $(C, S_C)$ is a log curve together with a non-constant
logarithmic map $\phi: (C,S_C)\to (\bar Y,\bar S)$, i.e $\phi^{-1}(\bar S)\subset S_C. $ Then $(C,S_C)$ is hyperbolic, i.e. $$\deg\Omega^1_C(\log S_C)>0.$$}\\[.2cm]
{\bf Proof:}
Applying Proposition 3.3.2.  on the pulled back family
$$\phi^*f:\phi^*X\to C,$$ which
is non-isotrivial and with the singular fibres $\phi^*\Delta$ over $S_C,$
we obtain a big sub sheaf
$$\mathcal A_C\hookrightarrow S^\ell\Omega^1_C(\log S_C),$$
which implies that $\Omega^1_{ C}(\log S_C)$ is big.\\[.4cm]
\subsection{ Augmented Line Bundles Attached to Families.}
Recall that Griffiths's augmented line bundle $\Lambda_{(E,\theta)}$ for a system of Hodge bunbdles
$$(E,\theta)=(\bigoplus_{p+q=n}E^{p,q},\bigoplus_{p+q=n}\theta^{p,q}) $$
arising from a polarized variation of Hodge structures over $\bar Y\setminus\bar S$
is defined as 
$$\Lambda_{(E,\theta)}:=(\det E^{n,0})^n\otimes(\det E^{n-1,1})^{n-1}\otimes\cdots\otimes (\det E^{0, n-1}).$$
It is ample over an open set of $\bar Y$ if the Torelli map is injective at a point.\\[.2cm] 
We introduce the Griffiths argumented line for
a family 
$$f:X\to Y$$
of varieties of semi-ample canonical line bundles as
$$\Lambda_f:=\left((\det R^nf_*T^n)^n\otimes(\det R^{n-1}f_* T^{n-1})^{n-1}\otimes\cdots\otimes (\det R^1f_*T)\right)^\vee,$$
where 
$$T^i:=\wedge^iT_{\bar X/\bar Y}(-\log\bar\Delta).$$ 
We remark that for a family of abelian varieties or Calabi-Yau manifolds the both
argumented line bundles coincide with each other.
\\[.2cm]
{\bf Proposition 4.2.1.}
{\sl Let $f:\bar X\to\bar Y$ be a family of $n$-folds of ample canonical line bundle and with maximal Var($f$). Then
	$\Lambda_f$ 
	is ample over an open set of $\bar Y.$}\\[.2cm]
{\bf Sketch the Proof.}\\
$\bullet$ Using the second comaprison map in Section 3.4 we find a
Kummer cover $\psi :\bar Y'\to Y$ such that there is a Higgs map
$$\rho:\psi^*(F,\tau)\to (E,\theta)\otimes A^{'-1}$$
over $\bar Y'$. In fact, $\rho$ is an inclusion follows from Prop 2.4 as
the canonical line bundle of the generic fibre of $f$ is ample.\\[.2cm]
$\bullet$ $\psi^*(F,\tau)$ has $n$ pieces of sub Higgs bundles of the form
$$\psi^*(G,\tau)_p=(\bigoplus_{i=p}^n\psi^*R^if_*T^i,\psi^*\tau),\quad 1\leq p\leq n.$$
And at least $(G,\tau)_n\not=0.$\\[.2cm]
$\bullet$
Assume the sub Higgs bundle
$$\psi^*(G,\tau)_p\subset (E,\theta)\otimes A^{'-1}$$
is non-zero. Then
$$\det\psi^*(G,\tau)_p\subset\wedge^r (E,\theta)\otimes A^{'-r}.$$
As $\det\psi^*(G,\tau)_p$ is a rank-1 Higgs bundle with the nilpotent Higgs field, and is the zero Higgs field.  It factors through
$$\det\psi^*(G_p)\subset\text{ker}(\wedge^r\theta)\otimes A^{'-r}.$$
This shows $\det\psi^*(G_p)^\vee$ is ample over an open set of $\bar Y'$
and $\det G_p^\vee $ is ample over an open set of $\bar Y$.\\
Putting all of them together we obtain
$$\left(\det G_n^n\otimes\det G_{n-1}^{n-1}\otimes\cdots\otimes\det G_1\right)^\vee=\left((\det R^nf_*T^n)^n\otimes\cdots\otimes (\det R^1f_*T)\right)^\vee=\Lambda_f
$$
is ample over an open set of $\bar Y.$
\\[.4cm]
\subsection{ Pseudo Effectivity and Viehweg Hyperbolicity.}
We recall the following conjecture asked by Viehweg\\[.2cm]
{\bf Conjectue 4.3.1.}(Viehweg) Let $f:\bar X\to\bar Y$ be a family of $n$-folds with the semi-ample canonical line bundle and with the quasi-finite classifying map. Then
$(\bar Y,\bar S)$ is of log general type.\\[.2cm]
{\bf Remark 4.3.2.} Conjecture 4.3.1. asks the so-called Viehweg hyperboilicity, i.e.
any log sub variety $(Z, S_Z)\subset (\bar Y,\bar S)$ is of log general type.  The conjecture was proven by Campana-Paun for the case of fibres with semi-ample canonical line bundles and by Popa-Schnell for the case of fibres of genral type or with good minimal model. \\[.2cm]
{\bf Theorem 4.3.3.}(Campana-Paun \cite{CP}, and Popa-Schnell \cite{PS}) {\sl 
Viehweg hyperbolicity holds true on moduli stacks of manifolds with good minimal model or of general type.}\\[.2cm]
Let's explain the crucial points in their proof.
Given a family 
$$f:\bar X\to \bar Y$$ 
 of varieties with good minimal model or of general type and with maximal variation by Kawamata, Viehweg and Koll\'ar one knows $f_*\omega_{\bar X/\bar Y}^\nu$ is big. Hence, by uning Griffiths-Yukawa coupling on the deformation Higgs bundle attached to $f$  one finds big sub sheaf $\mathcal A\subset S^\ell\Omega^1_{\bar Y}(\log\bar S)$, which 
induces a short exact sequence
$$0\to\mathcal A\to S^\ell\Omega^1_{\bar Y}(\log\bar S)\to Q\to 0,$$
The main ingredient in their proof is to show $\Omega^1_{\bar Y}(\log\bar S)$ has a suitable semi positivity, the so-called  {\sl pseudo effectvity} introduced by Demailly:\\[.2cm]
{\bf Definition-Proposition 4.3.4.} (Bouckksom-Demailly-Paun-Peternell)\\[.2cm]
{\sl A line bundle $L$ on a projective manifold $\bar Y$ is called pseudo-effective if
$L\cdot C\geq 0$ for all irreducible curves $C$ which move in a family covering $\bar Y.$}\\[.2cm]
{\bf Lemma 4.3.5.} (Demailly) {\sl Given two line bundles $L$ and $M$. If $L$ is pseudo effective and $M$ is big then $L\otimes M$ is big.}\\[.2cm]
{\bf Corollary 4.3.6.}{\sl The notion of pseudo-effective and the notion weakly positive for line bundles coincide with each other.}\\[.2cm]
{\bf Theorem 4.3.7.}(Campana-Paun) {\sl 
If some tensor power of $\Omega^1_{\bar Y}(\log\bar S)$ contains a sub sheaf of big
determinant. Then the determinant of every quotient sheaf of every tensor power
of $\Omega^1_{\bar Y}(\log\bar S)$ is pseudo-effective.}\\[.2cm]
Theorem 4.3.3. follows from Theorem 4.3.7. and Lemma 4.3.5.\\[.2cm]
{\bf Question 4.3.8.}(Viehweg-Zuo)
Let $f:V\to U$ be family of $n$-folds of semi-ample canonical line bundle and with
the quasi-finite classifying map into the moduli space, and let $\bar Y\supset U$ be a smooth projective compactification, we ask the following semi positivities, which are stronger than pseudo-effectivity.
\\[.2cm]
$\bullet$ For any morphism from a smooth projective curve 
$$\phi: \bar C\to \bar Y,\quad
\text{ such that}\quad \phi(\bar C)\cap U\not=\emptyset,$$  does any quotient bundle of
 $\phi^*\Omega^1_{\bar Y}(\log \bar S)$ 
 have non-negative degree?\\[.2cm]
$\bullet$ Does $\Omega^1_{\bar Y}(\log\bar S)$ be weakly positive?\\[.4cm]
\section{ \bf Picard Type Theorems and Relation to Hyperbolicity}
\subsection{ Classical Picard's Theorems.}
\'Emile Picard proved two important theorems in complex analysis:\\[.2cm]
{\bf Theorem 5.1.1.} (Little Picard Theorem)  {\sl 
There is no nonconstant holomorphic map 
$$\gamma: \mathbb{C}\to \mathbb{P}^1\setminus\{0,1,\infty\}.$$}\\[.2cm]
{\bf Theorem 5.1.2.} (Big Picard Theorem)
{\sl Any holomorphic map 
	$$\gamma :\mathbb{D}^*:=\{z\in\mathbb{C}\,;\,0<|z|<1\}\to \mathbb{P}^1\setminus\{0,1,\infty\}$$  can be extended as a holomorphic map from $\mathbb{D}$ to $\mathbb{P}^1$.}\\[.2cm]
Big Picard theorem is some times called Picard extension theorem.
\\[.4cm]
{\bf Notions of Complex Hyperbolicity in Higher Dimensions.}\\[.2cm]
$\bullet$ We say a complex space $X$ is {\sl Brody hyperbolic} if there is no nonconstant holomorphic map from $\mathbb{C}$ to $X$.\\[.2cm]
$\bullet$ We say a complex space $X$ is {\sl Kobayashi hyperbolic} if the Kobayashi pseudometric is nondegenerate.\\[.2cm]
If $X$ is proper then the both notions are equivalent. In general, Kobayashi hyperbolicity implies Brody hyperbolicity.\\[.2cm]
{\bf Theorem 5.1.3.}  (Kobayashi-Ochiai '70, A. Borel '72,.... ) {\sl 
A fine moduli space of polarized abelian varieties is Brody and Kobayashi hyperbolic, and the Picard extension theorem holds true by taking the Baily-Borel compactication.}\\[.2cm]
$\bullet$ Inspired by the above Picard extension theorem, Javanpeykar and Kucharczyk formulated the following notion:\\[.2cm]
A finite type scheme $X$ over $\mathbb{C}$ is called \emph{Borel hyperbolic} if, for every finite type reduced scheme $S$ over $\mathbb{C}$, any holomorphic map from $S$ to $X$ is algebraic.\\[.2cm]
{\bf Example:}\\
$\bullet$ $\exp:\mathbb C^*\to\mathbb C^*$
is holomorphic but not algebraic.\\[.2cm]
$\bullet$ Any holomorphic map $\gamma: S\to\mathbb C^*\setminus\{1\}$ is algebraic.
\\[.2cm]
We also recall the following theorem on the hyperbolicity of the period mapping:\\[.2cm]
{\bf Theorem 5.1.4.}  (Griffiths-Schmid) {\sl 
A quasi-projective manifold endowed with a locally injective period mapping is Brody and Kobayashi hyperbolic.}\\[.2cm]
The proof relies on the fundamental theorem due to Griffiths-Schmid
that the holomorphic sectional curvature of the Hodge metric is bounded above by a negative constant.
\\[.2cm]
{\bf Hyperbolicity on Moduli Stacks of Manifolds with Semi-Ample Canonical Line Bundle.}\\[.2cm]
Motived by the above theorems we study the various hyperbolicities on moduli stacks of manifolds of semi-ample canonical line bundles.\\[.2cm]
{\bf Theorem 5.1.5.} (not in the full generality yet) {\sl 
Let $U$ be a quasi-projective manifold with a smooth projective compactification $\bar Y\supset U.$ Let $f: V\to U$ is family of $n$-folds with semi-ample canonical sheaf $\omega,$  and assume the classifying map from $U$ to the moduli space is quasi finite. Then
\begin{itemize}
\item U is Brody hyperbolic
\item U is Kobayashi hyperbolic if $\omega$ is semi-ample and big.
\end{itemize}}.\\
The above theorem is due to several research groups.
\begin{itemize}
\item Viehweg-Zuo (2001) $U$ is Brody hyperbolic if $\omega$ is ample.\\
\item To-Yeung (2013) $U$ is Kobayashi hyperbolic
if $\omega$ is ample.\\
\item
Popa-Taji-Wu (2018) $U$ is Brody hyperbolic if $\omega$ is big.\\
\item Deng (2018) $U$ is Brody hyperbolic if $\omega$ is semi-ample and $U$ is Kobayashi hyperbolic if $\omega$ is semi-ample and big.

\end{itemize}.\\[.2cm]
{\bf Main idea:}
One constructs a complex finsler pseudo metric $h_F$ induced by the Hodge metric via Griffiths-Yukawa coupling of the Kodaira-Spencer map. $h_F$ has an extra  "strongly" negatively holomorphic sectional curvature arising from 
the semi-negativity of the kernel of Kodaira-Spencer map on a system of Hodge bundles plus an extral negativity coming from the dual of Kawamata-Viehweg's theorem on the positivity of the direct image of the relative pluri dualizing sheaf.\\
We like to point out, in the case of canonically polarized families To-Yeung \cite{TY} take the Weil-Petersson metric arising from K\"ahler-Einstein metric along the fibres, and run Griffiths-Yukawa coupling. The negativity of the kernel Kodaira-Spencer map comes from Siu's curvature computation for the Weil-Petersson metric.
\\[.4cm]
\subsection{\bf Picard Extension Theorem on Moduli Stacks}  The classical Picard extension theorem on $\mathbb C-\{0,1\}$ can be regarded as the Picard extension theorem on the fine moduli space of elliptic curves.  In this section we discusse Picard extension theorem for higher dimensional moduli stacks.\\[.2cm]
{\bf Theorem 5.2.1.}( Deng-Lu-Sun-Zuo \cite{DLSZ} Picard extension theorem on moduli stacks) {\sl 
Let $U$ be a quasi-projective manifold with a smooth projective compactification $\bar Y\supset U.$ Let $f: V\to U$ is family of $n$-folds with semi-ample canonical sheaf $\omega$ , and assume the classfying map from $U$ to the moduli space is quasi finite. Then Picard extension theorem holds true.}\\[.2cm]
{\bf Strategy}
\begin{itemize}
\item Constructing
a complex Finsler pseudo metric on $\bar Y$, which satisfies
a curvature inequality.
\item
A generalization of Lu's exension theorem
\end{itemize}.\\[.2cm]
First we recall the original Lu's extension theorem \cite{Lu}\\[.2cm]
{\bf Theorem 5.2.2.}
(Lu's extension theorem)  {\sl 
Let $\bar Y $ be a projective manifold and $\bar S $ a simple normal crossing divisor on $\bar Y $, and assume
some symmetric power of $\Omega^1_{\bar Y}(\log\bar S)$ contains an ample line bundle
$$A\hookrightarrow S^\ell\Omega^1_{\bar Y}(\log \bar S).$$
Let $\gamma:\,\mathbb{D}^*\to\bar Y\setminus\bar S $ be a holomorphic map such that the pulled back
$$\gamma^*A\hookrightarrow\gamma^*S^\ell\Omega^1_{\bar Y}(\log\bar S)\xrightarrow{d\gamma}S^\ell\Omega^1_{\mathbb D^*}$$
does not vanish,
then $\gamma$ extends to a holomorphic map from $\mathbb{D}$ to $\bar Y$.}\\[.2cm]
{\bf Picard extension theorem for moduli stacks of curves}.\\[.2cm]
Applying Lu's extension theorem we prove a baby case of the Picard extension theorem,
the Picard extension theorem for the moduli stack of projective curves. Of course, it follows also from Kobayashi-Ochiai-Borel's Picard extension theorem on the moduli space of abelian varieties via Jacobian map.\\[.2cm]
Let $$f: V\to\mathcal M_g$$ be the universal family of curves of genus $g$ over a fine moduli space $\mathcal M_g$ together with the Deligne-Mumford compactification
$$ f: X\to \bar {\mathcal M_g}.$$
The deformation Higgs bundle has the form
$$(\mathcal O\oplus R^1f_*T_{X/\bar{\mathcal M}_g}(-\log\bar\Delta),\tau)$$
with the isomorphic Kodaira-Spencer map $$\tau: T_{\bar{\mathcal M}_g}(-\log\bar S)\xrightarrow{\simeq}R^1f_*T_{X/\bar{\mathcal M}_g}(-\log\bar\Delta),$$
and the dual
$$\tau^\vee:  f_*\omega^2_{X/\bar{\mathcal M}_g}\xrightarrow{\simeq}\Omega^1_{\bar{\mathcal M}_g}(\log\bar S).$$
By 
Kawamata-Viehweg  $ f_*\omega^2_{X/\bar{\mathcal M}_g}$
is ample over $\mathcal M_g$, i.e.
there exists a map on $\bar{\mathcal M}_g$
$$\bigoplus A\to S^\ell f_*\omega^2_{X/\bar{\mathcal M}_g}$$
which is surjective over $\mathcal M_g.$\\[.2cm]
Given a non-constant holomorphic map $\gamma:\mathbb D^*\to\mathcal M_g$
one finds a direct summand $A$ such that the pulled back
$$\gamma^* A\to\gamma^* S^\ell\Omega^1_{\bar{\mathcal M}_g}(\log\bar S)\xrightarrow{d\gamma}S^\ell\Omega^1_{\mathbb D^*}$$
does not vanish.  Hence, Lu's extension theorem applies.\\[.2cm]
{\bf Discussion on Picard extension theorem for the general case}\\[.2cm]
Let $V\to U$ be a family of $n$-folds of semi-ample canonical line bundles and with the quasi-finite classifing map, and $\bar Y=U\cup\bar S$ a smooth compactification.\\[.2cm]
Recall { Comparison 3.2.},
there exists a Higgs map
$$\rho: (F,\tau)\to (E,\theta)\otimes A^{-1}$$
from the defomation Higgs bundle to
the system of Hodge bundles $(E,\theta)$ arising from the quasi-canonical extension of a PVHS from a smooth family $g: Z\setminus\Pi\to\bar Y\setminus (\bar S+\bar T).$\\[.2cm]
Let $$(G,\theta_G)=:\text{Im}(\rho)\subset (E,\theta)\otimes A^{-1},$$
then $G^{n,0}=\mathcal O_{\bar Y}$ and $\theta_G$ has logarithmic pole only along $\bar S.$
By {Torelli injectivity 3.3.1.} one finds an open sub variety $U_0\subset U$ such that
$$ T_{\bar Y}(-\log\bar S)|_{U_0}\otimes G^{n,0}|_{U_0}\xrightarrow{\theta^{n,0}_G}G^{n-1,1}|_{U_0}$$ is pointwisely 
injective on $U_0.$\\[.2cm]
Given a non-constant holomorphic map
$$\gamma: C:=\mathbb D^*\to U$$
we may assume $\gamma(C)$ is Zariski dense in $U$. Hence $\gamma(C)\cap U_0\not=\emptyset.$\\[.2cm]
Taking the pulled back via $\gamma$ we obtain those two Higgs bundles pulled back to $C$
$$\gamma^*(G,\theta_G)\subset\gamma^*(E,\theta)\otimes\gamma^*A^{-1}.$$
The pointwise Torell injectivity on $U_0$ allows us runing the maximal non-zero iteration of Kodaira-Spencer map along the map $\gamma$ and we obtain:\\[.2cm]
{\bf Proposition 5.2.3.}  {\sl There exists an positive integer $\ell_\gamma$ such that
$$T_C^{\ell_\gamma}\xrightarrow{(\gamma^*\theta_G)^{\ell_\gamma}\not=0}
\text{ker}(\gamma^*\theta_G^{n-\ell_\gamma,\ell_\gamma})\subset
\text{ker}(\gamma^*\theta^{n-\ell_\gamma,\ell_\gamma})\otimes\gamma^*A^{-1}.
$$}
\\[.2cm]
{\bf A complex Finsler pseudo metric $h_F$ on $\bar Y$.}\\[.2cm]
By runing the iteration of Kodaira-Spencer map on $\bar Y$ with the length $\ell_\gamma$ found in Proposition 5.2.3.
$$ S^{\ell_\gamma}T_{\bar Y}(-\log\bar S)\xrightarrow{\theta_G^{\ell_\gamma}}E^{n-\ell_\gamma,\ell_\gamma}\otimes A^{-1},$$
and taking the product of the Fubini-Study metric $g_A$ on $A$ and the Hodge metric $g_{hod}$ on $E^{n-\ell_\gamma,\ell_\gamma}$
we define the a Finsler pseudo metric by
$$ |v|^2_{h_F}:= |\theta_G^{\ell_\gamma}(v^{\otimes\ell_\gamma}) |^{2/m}_{g^{-1}_{A}\otimes g_{hod}}$$
for $ v\in T_{\bar{Y}}(-\log\bar{S}).$\\[.2cm]
{\bf Modification of $h_F$}.\\
Take an auxiliary function associated to the boundary divisor $\bar{S}=\sum\bar S_i$ and $\bar T=\sum\bar T_j,$ defined is as the product of the norm functions
$$l^\alpha=\prod_i (-\log |s_i|_{g_i}^2)^\alpha\cdot\prod_j(-\log |t_j|_{g_j}^2)^\alpha$$
where $s_i$ respectively $t_j$ is the canonical section
of the line bundle $L_i=\mathcal O_{\bar Y}(\bar S_i)$ repectively
$L_j=\mathcal O_{\bar Y}(\bar T_j)$ and $g_i$ respectively $g_j$ is a smooth Hermitian metric on $L_i$ respectively $L_j$.\\[.2cm]
We define a new singular hermitian metric $g_{\alpha}:= g_A\cdot l^{\alpha}$ on $A$ and hence the modified Finsler metric $h_{F,\alpha}$
$$
|v|^2_{h_{F,\alpha}}:= |\theta_G^{\ell_\gamma}(v^{\otimes\ell_\gamma}) |^{2/m}_{g^{-1}_\alpha\otimes g_{hod}}$$
for $ v\in T_{\bar{Y}}(-\log\bar{S}).$
\\[.2cm]
{\bf A curvature inequality for $h_{F,\alpha}$}.\\[.2cm]
$h_{F,\alpha}$ is bounded on $U$ follows from the work by Schmid and Cattani-Kaplan-Schmid on the asymptotic behavior of the Hodge metric near the boundary divisor under the quasi-canonical extension.\\[.2cm]
By Prop. 5.2.3. on the maximal non-zero iteration along $\gamma$, the pulled back
$\gamma^*h_{F,\alpha}$ via the map $$\gamma: C\to U$$ is a non-zero pseudo metric on $C. $  For $ v\in T_C$ it reads
as
$$ |v|^2_{\gamma^*h_{F,\alpha}}= |\gamma^*\theta_G^{\ell_\gamma}(v^{\otimes\ell_\gamma}) |^{2/m}_{\gamma^*g^{-1}_\alpha\otimes\gamma^*g_{hod}}$$ and 
such that
$$\gamma^*\theta_G^{\ell_\gamma}(v^{\otimes\ell_\gamma})
\in\text{ker}(\gamma^*\theta^{n-\ell_\gamma,\ell_\gamma})\otimes\gamma^*A^{-1}.
$$
i.e. the metric $\gamma^*h_{F,\alpha}$ is the product of the Hodge metric on $\text{ker}(\gamma^*\theta^{n-\ell_\gamma,\ell_\gamma})$ with the Fubini-Study metric on $\gamma^*A^{-1}.$\\
By Griffiths' curvature formular we show that the curvature form of $(\text{ker}(\gamma^*\theta^{n-\ell_\gamma,\ell_\gamma}),\gamma^*g_{hod})$
is semi-negative and obtain:\\[.2cm]
{\bf Theorem 5.2.4.} (A curvature inequality: Viehweg-Zuo, Deng-Lu-Sun-Zuo)  {\sl 
Let $$V\to U$$ be a family of smooth projective varieties of semi-ample canonical line bundles and with the quasi-finite classifying map. Let $\bar Y\supset U$ be a smooth compactification.
For a holomorphic map $$\gamma: C\to U$$ with the Zariski dense image then
the Finsler pseudo metric $h_{F,\alpha}$ constructed above for $\alpha>>0$ satisfies the following curvature inequality
$$
d d^c\log |\gamma'(z) |^2_{h_{F,\alpha}}\geq c\gamma^*\omega_{FS},
$$
where $c$ is a postive constant and $\omega_{FS}$ is the K\"ahler form of the Fubini-Study metric on $A$.}
\\[.2cm]
Griffiths and King \cite{GK} studied the higher dimensional generalization of value distribution theory, also known as Nevanlinna theory.
With it, they obtained a Nevanlinna-theoretic proof of Borel's theorem via negative curvature. The existence of a negativily curved complex Finsler pseudometric
in our setting makes this proof of Griffiths-King, suitably generalized, applicable.\\[.2cm]
{\bf Theorem 5.2.5.} (A generalization of Lu's extension theorem: Deng-Lu-Sun-Zuo) {\sl 
Let $(Z,D)$ be a log pair and $\gamma:\mathbb D^*\to Z\setminus D$
be a holomorphic map. Assume there exists a Finsler pseudo metric on $Z$ satisfying the above curvature inequality w.r.t. $\gamma$  stated in Theorem 5.2.4. Then $\gamma$ extends to a holomorphic map $\bar\gamma:\mathbb D\to Z$.}\\[.2cm]
{\bf Sketch the proof:}
We identify our $\mathbb{D}^*$ with the inverted punctured unit disk $\mathbb{D}^\circ:=\{z\in\mathbb{C};\text{$|z|\geq 1$}\}$ and set $\mathbb{D}_{r_0,r}:=\{r_0\leq |z| <r\}$ for a fixed
$r_0>1$. Denote by $\gamma:\,\mathbb{D}^\circ\to U$ the analytic map in question. Then we want to show that $\gamma$ extends over the point at infinity.\\[.2cm]
We consider the following Nevanlinna characteristic function for ${\gamma^*\omega_{FS}}$ defined by
$$
T_{\gamma^*\omega_{FS}}(r) :=\int^r_{r_0}\frac{d\rho}{\rho}\int_{\mathbb{D}_{r_0,\rho}}\gamma^*\omega_{FS}.$$
It is well-known,  in the following lemma that the asymptotic behavior of $T_{\gamma^*\omega_{FS}}(r)$ as $r\to\infty$ characterizes whether $\gamma$ can be extended over the $\infty$.\\[.2cm]
{\bf Lemma 5.2.6.}  (Well known, see Demailly, Noguchi-Winkelmann,..)
 {\sl $$T_{\gamma^*\omega_{FS}}(r)=O(\log r)$$
  if and only if $\gamma$ can be extended holomorphically over $\infty$}.\\[.2cm]
We show in our situation\\[.2cm]
{\bf Lemma 5.2.7.} $$T_{\gamma^*\omega_{FS}}(r)=O(\log r)$$
\\[.2cm]
{\bf Sketch the proof:} The above curvature inequality implies
$$ T_{\gamma^*\omega_{FS}}(r)\le\int^r_{r_0}\frac{d\rho}{\rho}\int_{\mathbb{D}_{r_0,\rho}}dd^c\log\lvert\gamma'(z)\lvert^2_{h_{F,\alpha}},
$$
and one also uses the Nevanlinna characteristic function of $(1,1)$-form of the pulled back Finsler metric, the holomorphic sectioan curvature is bounded above by a negative constant \\$\cdots\cdots$ and shows
$$\int^r_{r_0}\frac{d\rho}{\rho}\int_{\mathbb{D}_{r_0,\rho}}dd^c\log\lvert\gamma'(z)\lvert^2_{h_{F,\alpha}}=O(\log r).$$\\[.4cm]
\subsection{ \bf Picard Extension for Period Mapping.}
For a quas-projective manifold $U=\bar Y\setminus\bar S$ carries a locally injective Torelli mapping Deng
applied the above extension thoerem to the maximal non-zero iteration of cup product with Kodaira-Spencer map on the Griffiths' agument line bundle $\Lambda$ as the initial bundle in a suitable tensor product of the original system of Hodge bundles
$$S^\ell T_{\bar Y}(-\log\bar S)\xrightarrow{\theta^\ell}\text{ker}(\theta)\otimes\Lambda^{-1}.$$
The Hodge metric on the right side induces a complex Finsler   pseudo metric satisfying the required curvature inequality for Picard extension theorem.\\[.2cm]
{\bf Theorem 5.3.1.} (Deng) {\sl Assume $U=\bar Y\setminus\bar S$ carries a locally injective Torelli mapping, then
the Picard extension theorem holds true for $(\bar Y,\bar S).$}\\[.4cm]
\section{\bf Arakelov Inequality.}
\subsection{\bf Viehweg Line Bundle and Arakelov Inequality.}
Let $ M=M_h $ denote the coarse moduli space of canonically polarized manifolds and with a fixed Hilbert polynomial $h$. We take a good projective compactification and consider a log map $$\phi:(\bar Y,\bar S)
\to (\bar M,\bar S_M)$$ from a log curve.\\[.2cm]
Having such a negatively curved Finsler metric on $M$ in Section 5 the general principle Yau's form of Ahlfors-Schwarz Lemma suggests that there must be an inequailty between the logarithmic hyperbolic metric on $(Y, S) $ and some singular K\"ahler metric on $\bar M.$\\[.2cm]
Viehweg has constructed a series of classes of invertible sheaves $\lambda_\nu$ on $\bar M$, which are nef on $\bar M$ and ample on $M.$ They are natural in the sense that
if $$\phi: (\bar Y,\bar S)\to (\bar M,\bar S_M)$$ is induced by a semi-stable family $f:\bar X\to\bar Y$ over a   1-dimensional base  $\bar Y$  then $$\phi^*\lambda_\nu=\det f_*\omega_{\bar X/\bar Y}^\nu.$$
{\bf Theorem 6.1.1.} (Arakelov inequality for relative pluri canonical linear system, Viehweg-Zuo, Moeller-Viehweg-Zuo)  {\sl 
Let $f:\bar X\to\bar Y$ be a   birational non-isotrivial semit-stable family of $n$-folds over a smooth projective curve $\bar Y$ and with the bad reduction over $\bar S$. Then
$$\mu( f_*\omega_{\bar X/\bar Y}^\nu)
\leq{n\nu\over 2}\cdot\deg\Omega^1_{\bar Y}(\log\bar S),\quad\forall\nu\in\mathbb N,\,\text{with}\, f_*\omega^\nu_{\bar X/\bar Y}\not=0.$$
In general,
$$\mu(A)\leq{n\nu\over 2}\cdot\deg\Omega^1_{\bar Y}(\log\bar S),\quad\forall\,\text{non-zero sub bundle}\,\,
A\subset f_*\omega^\nu_{X/Y}.
$$}
Many people have made contributions to Arakelov inequalities:\\[.2cm]
$\bullet$ Arakelov-Parshin, Faltings and Deligne (sharp form):  Arakelov inequality for families of abelian varieties.\\[.2cm]
$\bullet $ Viehweg-Zuo: Arakelov euality holds for some special families of abelian varieties. They are precisely the families of abelian varieties over Shimura curves with Mumford-Tate group of Mumford-type.\\[.2cm]
$\bullet$ Tan \cite{Tan} and Liu \cite{Liu}:  Arakelov inequality holds strictly for families of curves of genus >1 
and $\nu=1$.\\[.2cm]
$\bullet$ Green-Griffiths-Kerr \cite{ggk}, Jost-Zuo \cite{JZ}, Peters \cite{Pet}, Viehweg-Zuo \cite{VZ06}and a very recent work by Biquard-Collier-Garcia-Prada-Toledo \cite{bcgt}
have studied Arakelov inequality for systems of Hodge bundles over curves.\\[.2cm]
$\bullet$ Viehwge-Zuo, Moeller-Viehweg-Zuo: Some results over higher dimensional bases, but are not very satisfied.\\[.2cm]
 We list the main steps in the proof for Theorem 6.1.1. without mentioning on technical difficulties.\\[.2cm]
$\bullet$ Given a sub bundle
$A\subset f_*\omega^\nu_{\bar X/\bar Y}$ of
$\rk A=r$, by taking
$\det A\subset (f_*\omega^\nu_{\bar X/\bar Y})^{\otimes r}$
and suitable power of self-fibre product of $f$
we may assume $A$ is a line bundle. We want to show
$$\deg A\leq{n\nu\over 2}\cdot\deg\Omega^1_{\bar Y}(\log\bar S).$$
$\bullet $ Assume $\deg A>0$, and after a base change $\bar Y'\to\bar Y$ we may assume $A=A'^\nu$ for a line bundle $A'$ over $\bar Y$.\\[.2cm]
$\bullet$ Take the $\nu$-th cyclic cover $\bar Z\to\bar X$
of the natural section $$s:\mathcal O_{\bar X}\to (\omega_{\bar X/\bar Y}\otimes f^* A^{'-1})^\nu.$$
$\bullet$
As in Section 3.2. we obtain a Higgs map
$$\rho: (F,\tau)\otimes A'\to (E,\theta),$$
where $ (E,\theta)$ comes from the system of Hodge bundles of a quasi-canonical extension a VHS of weight-$n$ from geometry origin over $\bar Y\setminus (\bar S+\bar T)$. After a semi-stable reduction we
may even assume it is a canonical extension.\\[.2cm]
$\bullet$ Noting that
$F^{n,0}\otimes A'=A'\subset E^{n,0}$ and for
$T:=T_{\bar Y}(-\log\bar S)$,
the sub line bundle $A'$ with $\theta$-action generats a sub Higgs bundle of the form
$$ (G,\theta_G)=({A'}^\text{sat}\oplus{\theta(T\otimes A')}^\text{sat}\oplus\cdots\oplus{\theta^\ell(T^\ell\otimes A')}^\text{sat},\theta)\subset (E,\theta),$$
where ${\theta^i(T^i\otimes A')}^\text{sat}$ is the saturation of $\theta^i(T^i\otimes A')\subset E.$\\[.2cm]
$\bullet $ Applying Simpson's polystability for the sub Higgs bundle $(G,\theta_G)\subset (E,\theta)$  we show
$$(\ell+1)\deg A'+{\ell(\ell+1)\over 2}\deg T\leq\deg G\leq 0,$$
and
$${\deg A\over\nu}=\deg A'\leq{\ell\over 2}\deg\Omega^1_{\bar Y}(\log\bar S)\leq
{n\over 2}\deg\Omega^1_{\bar Y}(\log\bar S).$$\\[.4cm]
\subsection{\bf Strict Arakelov Inequality of Families of Varieties of General Type.}
In contras, for families of varieties of general type we show
the Arakelov inequality always holds strictly for large $\nu:$\\[.2cm]
{\bf Theorem 6.2.1.}  (Lu-Yang-Zuo)  {\sl 
Let $f:\bar X\to\bar Y$ be a birational non-isotrivial semistable family of $n$-folds of general type
then
$$\mu (f_*\omega^\nu_{\bar X/\bar Y})
<{n\nu\over 2}\cdot\deg\Omega^1_{\bar Y}(\log\bar S)$$
$\forall\,\nu\in\mathbb N$ such that the linear system of
$\omega^\nu_{\bar X/\bar Y}$ defines a birational $\bar Y$-morphism $$\eta: X\to\mathbb P(A).$$
}\\[.2cm]
{\bf Sketch the proof for Theorem 6.2.1.}\\
$\bullet$ Assume there exists such an $\nu$ with the Arakelov equality
$$\mu_0:=\mu (f_*\omega^\nu_{\bar X/\bar Y})={n\nu\over 2}\cdot\deg\Omega^1_{\bar Y}(\log\bar S).$$
As for any sub bundle $A\subset f_*\omega^\nu_{\bar X/\bar Y}$ we have shown in Theorem 6.1.1. 
$$\mu A\leq{n\nu\over 2}\cdot\deg\Omega^1_{\bar Y}(\log\bar S)=\mu_0, $$
i.e. $f_*\omega^\nu_{\bar X/\bar Y}$ is a semi-stable vector bundle over $\bar Y$.\\[.2cm]
$\bullet$ We consider the $d$-th multiplication map
$$K_{m_d}\subset S^d(f_*\omega^\nu_{\bar X/\bar Y})\xrightarrow{m_d}I_{m_d}\subset f_*\omega_{\bar X/\bar Y}^{d\nu}$$
where $K_{m_d}$ is the kernel of the $d$-th multiplication map.\\[.2cm]
{\bf Observation:} $K_{m_d}$ is the sub space of all homogenuos plynomials of degree $d$ in the ideal defining the birational embedding of the fibres of $f.$
\\[.2cm]
Applying Arakelov inequality for the image of the $d$-multiplication
$I_{m_d}\subset f_*\omega_{\bar X/\bar Y}^{d\nu}$ we show
$\mu I_{m_d}\leq d\cdot\mu_0$.\\[.2cm]
On the other hand
$S^d(f_*\omega^\nu_{\bar X/\bar Y})$ is semi-stable of
slope $d\cdot\mu_0$ and $I_{m_d}$ is a quotient bundle we obtain $$\mu I_{m_d}= d\cdot\mu_0=\mu S^d(f_*\omega^\nu_{\bar X/\bar Y})=\mu K_{m_d}$$
and $K_{m_d}$ is a semi-stable sub bundle of $S^d(f_*\omega^\nu_{\bar X/\bar Y})$ of the same slope.\\[.2cm]
$\bullet$
After a base change of $\bar Y$ and twisting with a line bundle we may assume
$f_*\omega^\nu_{\bar X/\bar Y}$ is semi-stable of degree zero and
$S^d(f_*\omega^\nu_{\bar X/\bar Y})$ contains
$K_{m_d}$ as a semi-stabel sub bundle of degree zero.\\[.2cm]
{\bf Theorem 6.2.2.} (Simpson)  {\sl 
Let $\mathcal C_{dR}$ be the category of vertor bundles over $\bar Y$ with integrable connections and $\mathcal C_{Dol}$ be the category of semi-stable Higgs bundle of degree $0$. Then there exists an equivalent functor
$$\mathcal F:\mathcal C_{Dol}\rightarrow\mathcal C_{dR}.$$}
We just recall some properties about this functor.
Let $(E, 0)$ be a semistable Higgs bundle with the trivial Higgs field and of degree $0$ and Let $(E',0) $ be a sub Higgs bundle of $(E, 0)$ of degree $0.$\\[.2cm]
(1). The functor $\mathcal F$ preserves the tensor products. In particular it also preserves symmetric powers.\\[.2cm]
(2). The underlying bundle of the bundle $\mathcal F((E,0))$ with
the integrable connection is isomorphic to $E$. We call the connection to be canonical and denote it by $\nabla_{can}(E)$.\\[.2cm]
(3). The connection $\nabla_{can}(E)$ preserves $E'$ and $\nabla_{can}(E)\mid_{E'}=\nabla_{can}(E')$.\\[.2cm]
Consequently,\\
(4) For a semi-stable vector bundle $V$ of degree $0$. Then there exists an integrable connection $\nabla$ on $V$ such that for any $d\geq1$ and any sub vector bundle $K\subset S^d(V)$ of degree $0$, the connection $S^d(\nabla)$ on $S^d(V)$ preserves $K$.\\[.2cm]
$\bullet$ Applying (4) in our situation we find an integrable connection
$(f_*\omega^\nu_{\bar X/\bar Y},\nabla) $ such that
$S^d(\nabla)$ preserves $$K_{m_d}\subset S^d(f_*\omega^\nu_{\bar X/\bar Y})$$ $\forall d\in\mathbb N.$\\
i.e. for each point in $U$ we find an analytic open neighborhood
$\Delta\subset U$ and a flat base
$\mathbb V$ for the slutions of $(f_*\omega^\nu_{\bar X/\bar Y},\nabla)_{\Delta}$
and such that $K_{m_d}\subset S^d(f_*\omega^\nu_{\bar X/\bar Y})$
is spanned by a flat sub space $\mathbb K_{m_d}\subset S^d(\mathbb V).$\\[.2cm]
This means that we find a base of ${f_*\omega^\nu_{\bar X/\bar Y}}$ over ${\Delta}$ such that the coeffiecents of all polynomials in the ideals
defining the smooth fibres of $$f:\bar X\to\bar Y$$ under the birational embedding $$ X\to\mathbb P(f_*\omega^\nu_{\bar X/\bar Y})$$ are locally constant. A contradiction to $f$ is   birational non-isortrivial.\\[.4cm]
\section{\bf Geometry on Moduli Stacks of Varieties with Semi-Ample Canonical Line Bundles.}
\subsection{\bf Pointed Rigidity Problem.}
Given a fixed log pair $(\bar Y,\bar S)$ and let $M=M((\bar Y,\bar S), h)$ be the set of isomorphism classes of families $$f: V\to\bar Y\setminus\bar S$$ of $n$-folds of semi-ample canonical line bundle and with the fixed Hilbert polynomial $h.$\\[.2cm]
The generalized Shafarevich program over function fields asks
\begin{itemize}
\item
finiteness of $M.$
\end{itemize}
This problem reduces to the following two sub problems
\begin{itemize}
\item boundedness of $M$
\item rigidity of $M$.
\end{itemize}
The boundedness of $M$ has been shown by various people, for families of abelian varieties is due Faltings, for famiies over a 1-dim base is due to Bedulev-Viehweg and Viehweg--Zuo. Finally, Kov\'acs--Lieblich showed the boundedness for families over a higher dimensional base.\\
The rigidity problem of families of higher dimensional fibres is quite subtle as there exists nonrigid families of higher dimensional varieties simply by taking the product of two families of lower dimensional fibres.
But under certain additional assumptions we expect the rigidity property to hold true.
As an intial step,
Javanpeykar-Sun-Zuo recently showed that the so-called $N$-pointed
Shafarevich conjecture holds true.\\[.2cm]
{\bf Theorem 7.1.1.}  (N-pointed rigidity, Javanpeykar-Sun-Zuo) {\sl 
Let $$f\colon V\to U$$ be a family of $n$-folds of semi-ample canonical line bundle and with the quasi-finite classifying map. Let $C$ be a smooth quasi-projective curve with the smooth compactification $\bar C=C\cup\bar S$ containing points $\{c_1,\cdots, c_N\}\subset C $ with $$N\geq{n-1\over 2}\deg\Omega^1_{\bar C}(\log\bar S).$$ Then for any set of fixed points $\{u_1,\cdots,u_N\}\subset U$ the set of isomorphism classes of non-constant morphisms $$\phi\colon C\to U$$ with $\phi(c_1)=u_1,\cdots\phi(c_N)=u_N$ is finite.}\\[.2cm]
	{\bf Sketch the proof for Thoerem 7.1.1.}
A non-trivial deformation of a map with $N$-fixed points
$$\phi: (C,\{c_1,\cdots, c_N\})\to (U,\{u_1,\cdots u_N\})$$
will produce a sub line bundle of degree $N$ in the pulled back of the deformation Higgs bundle
$$L=\mathcal{O}_{\bar C}(c_1+...+c_N)\subset F^{n-1,1}_{\bar C}.$$
(Via a Kawamata base change) we get the comparison map
$$\rho: (F,\tau)_{\bar C}\to (E,\theta)_{\bar C}\otimes A^{-1}_{\bar C}.$$
The maximal non-zero iteration of Kodaira-Spencer map starting with the sub line bundle $L$ generates a sub Higgs bundle in $(E,\theta)_{\bar C}\otimes A^{-1}_{\bar C}$,
and Simpson's polystability shows the expected upper bound for $N.$\\[.2cm]
The $N$-pointed rigidity has found already applications in the {\sl Mordellicity} of $U$,  a notion introduced by S. Lang.\\[.2cm]
Very recently, Javanpeykar, Lu, Sun and Zuo are making new progress and we expect that the one-pointed Shafarevich conjecture holds true, it should follow from the existence of the complex Finsler metric with the negatively curved curvature bounded above by a negative (1,1)-form on $U$   discussed in subsection 5.2.
\\[.4cm]
\subsection{\bf Question on
Tensor Decomposition of Deformation Higgs Bundle over a Product Base}
We like to study non-rigid familes. As a motivation we
take a family $$f\colon X\to Y$$ of smooth surfaces as the product of
two families of curves over base curves
$$f: X=X_1\times X_2\xrightarrow{f_1,f_2}Y_1\times Y_2=Y.$$
Then the deformation Higgs bundle
of the product family $f$ is decomposed as a tensor product of deformation Higgs bundels of $f_1$ and $f_2$
$$ (F,\tau)=p_1^*(F,\tau)_{1}\otimes p_2*(F,\tau)_{2}.$$
Further more, It induces a decomposition on the second relative pluricanonical linear system
$$f_*\omega^2_{X/Y}=F_f^{0,2\vee}=p_1^*F^{0,1\vee}_{1}\otimes p_2^*F^{0,1\vee}_{2}$$
$$=p^*f_{1}(\omega^2_{X_1/Y_1})\otimes
q^*f_{2}(\omega^2_{X_2/Y_2}).$$
We see that the tensor product of the pluricanonical linear system has a specific geometric meaning, it reconstructs the splitting on the fibres of $f$.\\[.2cm]
Conversely, we like to study the structure of a non-rigid family.
Given a family  $$ f: X\to Y_1,$$ if it is non-rigid then it will be extended to a family
$$f: X\to Y_1\times Y_2$$
of the maximal variation. It it is still non-rigid, we continue to procedure,... and obtain a family over a product base
$$f: X\to Y_1\times ...\times Y_\ell=Y$$
birationally.\\
Let $f: X\to M$ be the universal family over the moduli stack of $n$-folds with semi-ample canonical sheaves, we consider the length $\ell_M$ of the Griffiths-Yukawa coupling on the universal deformation Higg
bundle
$$ S^{\ell_M} T_{\bar M}(-\log S)\xrightarrow{\tau^{\ell_M}\not=0}\text{ker}(\tau^{n-\ell_M,\ell_M}).$$
 We expect the length $\ell_M$ should look like the rank of a locally bounded symmetric domain.\\[.2cm]
{\bf Theorem 7.2.1.} (Product of sub varieties in moduli stacks, Viehweg-Zuo)\\[.1cm]  {\sl 
	\begin{itemize}
		\item If there exists a generically finite   rational map 
		$$ (Y_1, S_1)\times \cdots\times (Y_\ell, S_\ell)\to (\bar M, S)$$
		then $\ell \leq \ell_M$.\\
		\item
		If $\ell=\ell_M$, then
		$$\dim Y_i\leq \dim H^\ell (X_y,  T^{\ell}_{X_y}).$$
	\end{itemize}}.\\[.2cm]
The reuslt is far being optimal, we raise \\[.2cm]
{\bf Conjecture 7.1.2.}
		Put $T_{Y_i}(-\log S_i)=T_i, $ then   the tensor product of the first Kodaira-Spencer map along the   tangent bundles of the factors
		$$\bigoplus_{1\leq i_1<...<i_k\leq\ell}
		T_{i_1}\boxtimes...\boxtimes T_{i_k}\xrightarrow{\tau_{i_1}\otimes...\otimes\tau_{i_k}} (R^1f_*T_{X/Y}(-\log S))^{\otimes k}
		\xrightarrow{\cup}
		R^kf_*T^k_{X/Y}(-\log S)$$
		is injective. In particular,
		$$\sum_{1\leq i_1<...<i_k\leq\ell}\dim Y_{i_1}\cdots\dim Y_{i_k}\leq
		\dim H^k(X_y, T^k_{X_y}).$$\\[.2cm]
We even ask for the following stronger form.\\[.2cm]
{\bf Conjecture 7.2.3.} There exist graded logarithmic Higgs bundles $ (F,\tau)_i$ over $(Y_i, S_i)$
such that
$$(F,\tau)_Y\supset (F,\tau)_1\boxtimes\cdots\boxtimes (F,\tau)_\ell.$$
\\[.2cm]
Conjecture 7.2.3. holds true for families of abelian varieties or Calabi-Yau varieties. It
follows from Deligne's tensor product decomposition of irreducible PVHS over product bases.\\[.2cm]
{\bf Proposition 7.2.4.} (Viehweg-Zuo \cite{VZ-4}) {\sl Conjecture 7.2.3. holds true for normalized famlies of hypersurfaces in $\mathbb P^n$ of non-negative Kodaira-dimension, and 
Conjecture 7.1.1. holds true for this case.}\\[.4cm]
{\bf Sketch the proof} \\
$\bullet$   Consider a normalized family of hypersurfaces in $\mathbb{P}^n$
$$ X\hookrightarrow\mathbb{P}^n\times Y=\mathbb P^n_Y$$
of degree $d\geq (n+1),$
then the line bundle $\mathcal O_{\mathbb P^n_Y}(X)$ is $d$-divisible.\\[.2cm]
$\bullet$
Take the  $d$-th cyclic cover $$Z\to\mathbb P^n_Y$$ ramified along $X$ we obtain the family $$g: Z\to Y\to Y$$
 parameterizing
the $d$-cyclic covers of $\mathbb P^n_y$ ramified along $X_y$ together with an $\mathbb{Z}/d\mathbb{Z}$-action.\\[.2cm]
$\bullet$ Taking PVHS of the local system $\mathbb W=R^{n}g_*(\mathbb Z_{Z_0})$  of the middle cohomology
we choose the $i$-th eigen-space
$\mathbb W_{i}$ for $i=d-(n+1)$.  Then $\mathbb W_i$ is a PVHS of Calabi-Yau type.\\[.2cm]
$\bullet$ Take the grading of the Hodge filtration on Deligne's quasi-canonical extension 
of $\mathbb{W}_{i}$  
$$(G,\theta)_i =(\bigoplus_{p+q=n}G_i^{p,q},\bigoplus_{p+q=n}\theta_i^{p,q})$$
over $\bar Y$ with the logarithmic pole along $\bar S.$  Here 
$\mathbb W$ is called of Calabi-Yau type, if  
$G^{n,0}\simeq \mathcal O$ and
$$\theta^{n,0}:G^{n,0}\otimes T_{\bar Y}(\log \bar S)\to G^{n-1,1}$$
is an inclusion.
The residue map around the divisor $X\subset \mathbb P^n_Y$ induces an inclusion of Higgs bundles 
$$(\bigoplus_{p+q=n, q\not=n}G_i^{p,q},\bigoplus_{p+q=n, q\not =n}\theta_i^{p,q})\xrightarrow{\text{Res}_X} (F,\tau)$$
$\bullet$ For $\mathbb W_i$ over $U_1\times\cdots\times U_\ell$ and applying   Deligne's tensor product decomposition for the irredusible 
summand $\mathbb W'_i\subset \mathbb W_i$ containing $G^{n,0}_i\simeq \mathcal O$ we obtain
the decomposition on $(F,\tau)$  conjectured in 7.1.3. via $\text{Res}_X$.\\[.4cm]
\subsection{ Question on Characterization of Fibres of Non-Rigid Families}.\\[.2cm]
	{\bf Conjecture 7.3.1.}  A non-rigid family $f: X\to Y_1\times Y_2$ of $n$-folds of ample canonical sheaves
and with the maximal variation is decomposed as a product of two families of lower dimensional fibres over base $Y_1$ and $Y_2$ (up to
a domained map).\\[.2cm]
There could be some hopes for families of hypersurfaces in projective spaces. For simplicity we consider 
hypersurfaces in $\mathbb P^3$  of dgree $d \geq 4.$\\[.2cm]
{\bf Proposition 7.3.2.} (Viehweg-Zuo)  {\sl Given two homogneous polynomials of two variables
	$f_d(x_1,x_2)$ and $ g_d(y_1,y_2)$ of degree $d$ without multiple roots.
	 Then the hypersurface $X_{f_d, g_d}\subset \mathbb P^3$ defined by 
	$$ f(x_1,x_2, y_1, y_2)=f_d(x_1, x_2)+g_d(y_1, y_2)$$ 
	is smooth, and is rationally domained by the product of two plane curves
	$$C_1 : x_3^d=f_d(x_1,x_2), \quad C_2: y_3^d=g_d(y_1,y_2).$$}\\[.2cm]
Summer up: the family  $\{X_{f_d, g_d}\}$  is parameterized by the pair of the roots of $f_d$ and the roots of $g_d$. In this way, we construct a non-rigid family of hypersurfaces in $\mathbb P^3$
$$ f: X\to \Sigma_d\times \Sigma_d$$
where  $\Sigma_d$ is the moduli space of
$d$ distingushed points in $\mathbb P^1.$ Moreover,
the family $f$ is rational domained by the product of two copies of these families 
of plane curves as  $d$-cyclic covers of $\mathbb P^1$.
\\[.2cm]
{\bf Conjecture 7.3.3.}  Any non-rigid families of hypersurfaces in $\mathbb P^3$ has the above form upto a projective transformation.\\[.2cm]
In fact,  in the proof for Proposition 7.1.4. sketched above Viehweg-Zuo have found a close relationship between the relative graded Jacobian ring and a system of Hodge bundles by showing that the graded pieces of the Jocobian ring of a hypersurface $X\subset \mathbb P^n$  coincide with the eigenspaces of the $\mathbb Z/d\mathbb Z$-action on the middle-dimensional Hodge cohomology of the cyclic cover 
of $\mathbb P^n$ ramified over $X$ as modules together with the Kodaira-Spencer multiplication. Using Deligne's result on tensor decomposition of a polarized variation of the middle cohomology of cyclic covers of hypersurfaces in the projective space $\mathbb P^n$  leads to a tensor product decompostion of the reletive Jacobian ring as modules together with the Kodaira-Spencer multiplication. 
It is well-known that the Jacobian ring of a smooth hypersurface determines the isomorphism type of the hypersurface itself. We hope that the above decomposition of the relative Jacobian ring (as a module) as well as the decomposition of the deformation Higgs bundle will lead to some significant geometric consequences on the fibers, as in Conjecture 7.3.1. We pose the following problem which gives a criterion for rigidity.\\[.15cm]
{\bf Problem 7.3.4.}  Show that any family $$f\colon V\to U$$  of maximal variation and containing a fiber with the generically ample cotangent sheaf is   rigid.\\[.15cm]

\section{\bf Riemann-Finsler Metric on Moduli Spaces of Varieties with Semi-Ample Canonical Line Bundles and Topological Hyperbolicity}.\\
	One notes that a punctured Riemann surface $U$ is {hyperbolic} if and only if 
	$\pi_1(U)$ is {infinite and nonabelian}.  J. Milnor introduced a growth function to messure how non-abelian an infinite group could be.\\[.2cm]
	{\bf Definition 8.1.1} (Milnor  \cite{MIL} )  A {growth function} $\ell$ associated to a finitely generated group $G$ is defined as follows: for each positive integer $s$ let $\ell(s)$ be the number of distinct group elements which can be expressed as words of length $\leq s$ with a fixed choice of generators and their inverse.\\[.2cm]
	{\bf Theorem 8.1.2.} (Milnor \cite{MIL})  {\sl 
		The fundamental group a compact Riemannian manifold $M$ with all Riemannian sectional curvatures less than zero has exponential growth.}\\[.2cm]
The proof relies on G\"unther's volume comparison theorem on the exponential growth of the volume of the geodesic ball on the universal cover $\tilde M$. 
	\\[.2cm]
	Motived by Milnor's theorem we introduce the notion
\\[.2cm]
{\bf Definition 8.1.3.} 
A quasi-projective manifold $U$ is called topologically hyperbolic if
for any positively dimensional subvariety $ Z\subset U$ the image of the natural map $$\pi_1(Z)\to \pi_1(U)$$ grows exponentally.\\[.2cm]
	{\bf Conjecture 8.1.3.}
	Let $U$ be a base parameterizing manifolds with semi-ample canonical line bundle and with the quas-finite classifying map. Then $U$ is topologically hyperbolic.\\[.2cm]
	If $U$ carries an injective Torelli map, then $\pi_1(U)$ is non-abelian and infinite. One should be able to show it grows exponentially.
	The theorem below is an evidence for this conjecture without the assumption on the existence of an injective Torelli map.   \\[.2cm]
	{\bf Theorem 8.1.4.} (Lu-Sun-Zuo) {\sl The fundamental group of a fine moduli space of projective surfaces
		of Kodaira-dimension 1 and without multiple fibres in the Iitaka (elliptic) fibrations
		grows exponentially.}\\[.2cm]
	We propose a stratege to Conjecture 8.1.3. Similar to the approach in proving the complex hyperbolicity on $U$,  we hope to construct a Riemann-Finsler metric on $U$  via the iteration of Kodaira-Spencer map, whose curvature has certain negativity.
	\begin{itemize}
		\item {Wu-Xin} \cite{WX} generalized Milnor's theorem to the Finsler setting by showing that the G\"unther's volume comparison theorem \cite{Gue}  holds for Riemann-Finsler metrics with negative flag curvatures  (an analogue of the Riemannian sectional curvature in Finsler geometry) introduced by Chern and Bao-Chern-Shen \cite{BCS}
\item In our situation,  by taking a finite sum of the complex Finsler metrices
		from finitely many different cyclic covers we obtain a non-degenerated complex Finsler metric on $U$, which naturally induces a non-degenerated Riemann-Finsler metric $ds^2_{\mathbb{R}\text{Fin}}$ on $U$.
		\item We are aware that in general the Riemannian curvature decreasing principle does not hold true for real sub-manifolds. Very recently together with Lu and Sun we observed the fact that pluriharmonicity of the composition of the horizontal period map   with the projection to the symmetric space of non-compact type implies that decreasing Riemannian curvature still holds true in a weak form.
		\item We expect this weak form of the negative sectional curvatures can be used to show G\"unter's volume comparison inequality and get the solution of Conjecture 8.1.3. by Milnor's original argument.
	\end{itemize}

\end{document}